\documentclass[final]{siamltex}
\usepackage{amsmath}
\usepackage{amssymb}
\usepackage{amscd}
\usepackage{color}
\usepackage{amsfonts}
\usepackage{graphicx}
\usepackage{todonotes}
\usepackage{bm}
\usepackage{tikz}
\usepackage{pgfplots}
\usepackage{multirow}
\newcommand\x{\ensuremath{\bm{x}}}
\newcommand\z{\ensuremath{\bm{z}}}
\newcommand\bb{\ensuremath{\bm{b}}}

\author{
Siegfried Cools$^*$
\and Bram Reps$^*$
\and Wim Vanroose\thanks{Applied Mathematics Research Group, Department of Mathematics and Computer Science, University of Antwerp, Middelheimlaan 1, 2020 Antwerp, Belgium}
}

\title{An efficient multigrid calculation of the Far Field map for Helmholtz and Schr\"odinger equations}

\begin{document}
\maketitle

\begin{abstract}
  In this paper we present a new highly efficient calculation method
  for the far field amplitude pattern that arises from scattering
  problems governed by the $d$-dimensional Helmholtz equation
  and, by extension, Schr\"odinger's equation. The new technique is based 
  upon a reformulation of the classical real-valued Green's function 
  integral for the far field amplitude to an equivalent integral over a complex 
  domain. It is shown that the scattered wave, which is essential for the 
  calculation of the far field integral, can be computed very efficiently 
  along this complex contour (or manifold, in multiple dimensions). Using the iterative multigrid method 
  as a solver for the discretized damped scattered wave system, the proposed 
  approach results in a fast and scalable calculation method for the far 
  field map. The complex contour method is successfully validated on 
  Helmholtz and Schr\"odinger model problems in two and three spatial 
  dimensions, and multigrid convergence results are provided to substantiate 
  the wavenumber scalability and overall performance of the method.
\end{abstract}



\begin{DOI}
\end{DOI}

\pagestyle{myheadings}
\thispagestyle{plain}
\markboth{S. COOLS, B. REPS AND W. VANROOSE}{MULTIGRID CALCULATION OF THE FAR FIELD MAP}

\section{Introduction} \label{sec:Introduction} Scattering problems
are of key importance in many areas of science and engineering since
they carry information about an object of interest over large
distances, remote from the given target. Consequently, ever since
their original statement a variety of applications of scattering
problems have arisen in many different scientific subdomains.  In
chemistry and quantum physics, for example, virtually all knowledge
about the inner workings of a molecule has been obtained through
scattering experiments \cite{newton2002scattering}.  Similarly, in
many real-life electromagnetic or acoustic scattering problems
information about a far away object is obtained through radar or sonar
\cite{colton1998inverse}, intrinsically requiring the solution of 2D
or 3D wave equations. 

The near- and far field of a scattering problem are specific regions 
in space identified by the distance from the scattering object.
In these regions the scattering solution is subject to certain physical 
and mathematical properties. In the far field, the solution behaves as a 
spherical outgoing wave with an amplitude that decays as $1/r$, where 
$r$ is the distance to the scattering object. In this regime, the 
(asymptotic) solution can be written as a product of an angular part 
and a radial part. The near field is the region just outside the 
scattering object, where the solution has not yet reached its asymptotic 
form. In electromagnetic scattering, for example, the near field is 
dominated by dipole terms. Predicting the near- and far field solution 
is important for many present-day industrial applications (e.g.~radar), 
but also plays a key role in physical measurement systems such 
as tomography, near field microscopy and MRI.

However, far field maps are not limited to 2D or 3D scattering problems. 
New state-of-the-art experimental techniques in physics measure the full
4$\pi$ angular dependency of multiple fragments escaping from a
molecular reaction \cite{moshammer19964pi}. Through these experiments,
the reaction rates involving multiple fragments can be detected in
coincidence. Many experiments are being planned at facilities such as
e.g.\ the DESY Free-electron laser (FLASH) in Hamburg or the Linac
Coherent Light Source (LCLS) in Stanford. The accurate prediction of
the corresponding amplitudes starting from first principles requires
the use of efficient numerical methods to solve the high-dimensional
scattering problems, which can scale up to 6D or 9D in this context.
Indeed, after discretization one generally obtains a large, sparse and
indefinite system of equations in the unknown scattered wave. Direct
solution of this system is usually prohibited due to the massive size
of the problem in higher spatial dimensions. Preconditioned Krylov subspace
methods are able to solve some symmetric positive definite systems in
only $\mathcal{O}(n)$ iterations, where $n$ is the number of unknowns
in the system \cite{zubair2007multigrid,van2003iterative}. However,
scattering problems are often described by indefinite Helmholtz equations, 
which are generally hard to solve using iterative methods.  

This paper focuses on calculating the far field map resulting from
Helmholtz and Schr\"odinger type scattering problems \cite{baertschy2001electron}, 
which yields a $360^\circ$ representation of the scattered wave amplitude 
at large distances from the object of interest. The calculation of the far 
field map can typically be considered a two-step process. First, a Helmholtz 
problem with absorbing boundary conditions is solved on a finite numerical box
covering the object of interest. This step generally is the main computational 
bottleneck in far field map computation, since discretization of a Helmholtz 
equation leads to a highly indefinite linear system that is notoriously hard to
solve using the current generation of iterative methods. In particular, the
highly efficient iterative multigrid method 
\cite{brandt1977multi,briggs2000multigrid,trottenberg2001multigrid}
is known to break down when applied to these type of indefinite Helmholtz problems. 
The observed instability of both the coarse grid correction and relaxation scheme 
is due to close-to-zero eigenvalues of the discretized operator on some intermediate 
multigrid levels \cite{elman2002multigrid,ernst2010difficult}. In the second step a 
volume integral over the Green's function and the numerical solution is calculated 
to obtain the angular-dependent far field amplitude map. The two-step 
strategy was successfully applied to calculate impact ionization in
hydrogen \cite{rescigno1999collisional} and double photo-ionization in
molecules \cite{vanroose2006double} described by the Schr\"odinger equation, 
which in this case translates into a 6-dimensional Helmholtz problem.
However, computational overhead due to repeatedly solving the Helmholtz systems is significant, 
and supercomputing infrastructure was previously required to perform this type of calculations.

In this paper we propose a new method for the calculation of the far field
map that aims at bypassing the computational bottleneck in solving the Helmholtz equation. 
Using basic complex analysis, the method reformulates the far field integral over the Green's function on a
complex contour. The advantage of this reformulation is that the far field integral now requires the solution 
of the Helmholtz equation along a complex contour, which corresponds to a damped equation, instead of requiring the real-valued scattered wave solution.
The problem of solving a damped Helmholtz equation is well-known in the literature around Helmholtz preconditioners to be significantly easier than its real-valued counterpart.

Indeed, over the past decade significant research has been performed on the
construction of good preconditioners for Helmholtz problems. These results prove to be 
valuable in the context of this paper, albeit not in a preconditioning setting.
Recent work related to this topic includes the wave-ray approach \cite{brandt1997wave}, the
idea of separation of variables \cite{plessix2003separation},
algebraic multilevel methods \cite{bollhoefer2009algebraic}, multigrid
deflation \cite{sheikh2013convergence} and a transformation of the
Helmholtz equation into an advection-diffusion-reaction problem
\cite{haber2011fast}. Moreover, in 2004 the Complex Shifted Laplacian 
(CSL) was proposed by Erlangga, Vuik and Oosterlee
\cite{erlangga2006novel} as an effective preconditioner for Helmholtz problems. 
The key idea behind this technique is to formulate a
perturbed Helmholtz problem that includes a complex-valued wavenumber.  
Given a sufficiently large complex shift, this implies a
damping in the problem, thus making the perturbed problem solvable using multigrid in
contrast to the original Helmholtz problem with real-valued wavenumbers.
The concept of CSL has been further generalized in a variety of papers among which
\cite{erlangga2008multilevel,osei2010preconditioning}.
Recently a variation on the Complex Shifted Laplacian scheme by the
name of Complex Stretched Grid (CSG) was proposed in
\cite{reps2010indefinite}, introducing a 
complex-valued grid distance instead of a complex-valued wavenumber in the 
perturbed system. It was furthermore shown that the CSG system
is generally equivalent to the CSL scheme, and 
thus allows for a fast and scalable solution using classical multigrid methods.

These observations prove particularly useful in the context of far field map 
computation, since the damped Helmholtz equation appears naturally in the 
reformulation of the far field map integral proposed in this work.
The level of damping is governed by the size of the complex shift (CSL) or rotation 
angle (CSG), which is well-known to be bounded from below by the requirement of a stable 
multigrid solver, see \cite{cools2012local,erlangga2006novel,vangijzen2007spectral}. 
On the other hand, choosing the shift too large may negatively 
impact the accuracy of the integral quadrature, hence imposing an upper 
bound on the damping parameter.

To validate our approach, the method is successfully illustrated on both 2D and 3D Helmholtz 
and Schr\"odinger equations for a variety of discretization levels.
The absorbing boundary conditions used in this paper are based 
on the principle of Exterior Complex Scaling (ECS) that was introduced in the 1970's
\cite{aguilar1971class,balslev1971spectral,simon1979definition}, and is nowadays
frequently used in scattering applications. This method is equivalent to a
complex stretching implementation of Perfectly Matched Layers (PML) \cite{berenger1994perfectly,Chew20073d}.

The outline of the article is the following. In Section \ref{sec:Helmholtz} we introduce 
the notation and terminology that will be used throughout the text. Additionally, we
illustrate the classical calculation of the far field map for
Helmholtz type scattering problems. Section \ref{sec:Far} contains the main theoretical insights presented in this work. 
Here we introduce an alternative way of calculating the far field mapping
based upon a reformulation of the integral over a complex contour, for
which the corresponding Helmholtz system is very efficiently solved
iteratively. The new technique is validated in Section \ref{sec:Numerical}, and convergence 
results are shown for a variety of Helmholtz type model
scattering problems in both two and three spatial dimensions. We demonstrate
that the method allows for a very fast and scalable far field map calculation. 
In Section \ref{sec:Schrodinger}, the method is extensively tested on several semi-realistic Schr\"odinger 
type model problems in 2D and 3D respectively. Benchmark problems 
include quantum-mechanical model problem applications for which single, double and/or triple ionization occur. 
Finally, along with a discussion on the topic, conclusions are drawn in Section \ref{sec:Conclusions}.

%
%
%
%

\section{Helmholtz equation and far field map} \label{sec:Helmholtz}
In this section we introduce the general notation used throughout the text and we illustrate the classical derivation of the far field scattered wave solution and the calculation of its amplitude from a general Helmholtz type scattering problem, see \cite{colton1998inverse,kirsch1996introduction}. The theory presented in this section provides the foundation for the work in \cite{mccurdy2004solving}, where the below results were applied to a specific quantum-mechanical problem (cf.~Section \ref{sec:Schrodinger}).

%
%

\subsection{Notation and basic concepts} \label{subsec:notation}

The Helmholtz equation is a simple mathematical representation of the physics behind a wave scattering at a certain object. Let the object be defined on a compact support area $O$ located within a domain $\Omega \subset \mathbb{R}^d$. The equation is given by
\begin{equation} \label{eq:helmholtz_rhs}
 \left( -\Delta - k^2(\x)\right) u(\x) = f(\x), \qquad \x \in \mathbb{R}^d,
\end{equation}
with dimension $d \geq 1$, where $\Delta$ is the Laplace operator, $f$ designates the right hand side or source term, and $k$ is the (spatially dependent) wavenumber, representing the material properties inside the object of interest. Indeed, the wavenumber function $k$ is defined as
\begin{equation}
  k(\x) = \begin{cases}
  	k_1(\x), 		 \quad \textrm{for~} \x \in O,\\
  	k_0,   \hspace{0.16cm} \quad \quad \textrm{for~} \x \in \mathbb{R}^d \setminus O,
  \end{cases}
\end{equation}
where $k_0 \in \mathbb{R}$ is a scalar constant denoting the wavenumber 
outside the object of interest. The scattered wave solution is
given by the unknown function $u$. Throughout the text we will use the
following convenient notation
\begin{equation}
  \chi(\x) :=  \frac{k^2(\x) - k_0^2}{k_0^2},
\end{equation}
such that $k^2(\x)= k_0^2\left(1+ \chi(\x)\right)$. Note that the
function $\chi$ is trivially zero outside the object of interest $O$
where the space-dependent wavenumber $k(\x)$ is reduced to $k_0$. 
Defining the incoming wave as $u_{in}(\x) = e^{i k_0 \eta \cdot \x}$, 
where $\eta$ is a unit vector that defines the direction, 
the right-hand side is typically given by 
\begin{equation}
	f(\x) = k_0^2\chi(\x)u_{in}(\x). 
\end{equation}
The above expression follows directly from the fact that the total wave 
$u_{tot} = u + u_{in}$ satisfies the homogeneous Helmholtz equation 
$(-\Delta-k(\x)) \, u_{tot}(\x) = 0$, and the incoming wave trivially 
satisfies $(-\Delta-k_0^2) \, u_{in} = 0$. 
Reformulating (\ref{eq:helmholtz_rhs}), we obtain
\begin{equation} \label{eq:helmholtz_in}
  \left(-\Delta - k^2(\x)\right) u(\x) = k_0^2 \chi(\x)u_{in}(\x), \qquad \x \in \mathbb{R}^d.
\end{equation}
This equation is typically formulated on a bounded domain $\Omega$ with outgoing wave
boundary conditions on $\partial\Omega$, and can in principle be solved on a numerical box 
(i.e.\ a discretized subset of $\Omega^N \subset \Omega$) covering the support of $\chi$, with
absorbing boundary conditions along all edges. Let us assume that the numerical solution satisfying 
(\ref{eq:helmholtz_in}) has been calculated and is denoted by $u^N$.

\subsection{Classical derivation of the far field map} \label{subsec:classical}

In order to calculate the far field scattered wave pattern, the above equation is reorganized as
\begin{equation} \label{eq:helmholtz_reorganized}
  \left(-\Delta - k_0^2\right) \, u(\x) = k_0^2 \chi(\x)\left(u_{in}(\x)+u(\x)\right), \qquad \x \in \mathbb{R}^d.
\end{equation}
Note that we can replace the function $u(\x)$ in the right hand side of this equation with
the numerical solution $u^N(\x)$ obtained from equation \eqref{eq:helmholtz_in}.
In doing so, the above equation becomes an inhomogeneous Helmholtz equation with
constant wavenumber
\begin{equation} \label{eq:HHR}
  \left(-\Delta - k_0^2\right) \, u(\x) = g(\x),  \qquad \x \in \mathbb{R}^d,
\end{equation}
where the short notation $g(\x):=k_0^2 \chi(\x) (u_{in}(\x)+u^N(\x))$
is introduced for readability and notational convenience. It holds that $g(\x) = 0$ for $\x \in \mathbb{R}^d \setminus O$.
The above equation can easily be solved analytically using the Helmholtz Green's function
$G(\x, \x^\prime)$, i.e.
\begin{equation} \label{eq:scint}
	u(\x) = \int_{\mathbb{R}^d} G(\x, \x^\prime) \, g(\x^\prime) \, d\x' ,  \qquad \x \in \mathbb{R}^d.
\end{equation}
Since the function $g$ is only non-zero inside the numerical box
that was used to solve equation (\ref{eq:helmholtz_in}), the above
integral over $\mathbb{R}^d$ can be replaced by a finite integral over $\Omega$
\begin{equation}  \label{eq:first_int}
  u(\x) = \int_\Omega G(\x,\x^\prime) k_0^2 \chi(\x^\prime) \left(u_{in}(\x^\prime) + u^N(\x^\prime) \right) \, d\x^\prime, \qquad \x \in \mathbb{R}^d.
\end{equation}
In practice, this expression allows us to calculate the scattered wave solution $u$ in any point 
$\x \in \mathbb{R}^d\setminus \Omega^N$ outside the numerical box, using only the information $\x \in \Omega^N$ inside the numerical box.

Given the integral expression (\ref{eq:first_int}), the asymptotic form of the Green's function can be used to compute the far field map of the scattered wave $u$. 
In the following, this will be illustrated for a 2D model example where the Green's function is given explicitly by
\begin{equation}
  G(\x,\x^\prime) = \frac{i}{4}H_0^{(1)}(k_0|\x-\x^\prime|), \qquad \x, \x^\prime \in \mathbb{R}^d.
\end{equation}
where $i$ represents the imaginary unit and $H_0^{(1)}$ is the 0-th order Hankel function of the first kind. An analogous derivation can be performed in 3D, where we mention for completeness that the Green's function is given by
\begin{equation}
  G(\x,\x^\prime) = \frac{e^{i k_0 |\x - \x^\prime|}}{4\pi|\x - \x^\prime|}, \qquad \x, \x^\prime \in \mathbb{R}^d.
\end{equation}

To calculate the angular dependence of the far field map, the
direction of the unit vector $\boldsymbol\alpha$ is introduced that is
in 2D defined by a single angle $\alpha$ with the positive horizontal
axis, i.e.\ $\boldsymbol\alpha = (\cos \alpha, \sin
\alpha)^T$. Rewriting the spatial coordinates $\x$ in polar
coordinates as $\x = \left(\rho \cos\alpha, \rho\sin\alpha\right)^T$,
the asymptotic form of the Green's function for $|\x| \gg 1$ $(\rho
\to \infty)$ is given by
\begin{align}
  \frac{i}{4} H_0^{(1)}(k_0|\x-\x^\prime|) &= \frac{i}{4} \sqrt{\frac{2}{\pi}}e^{-i\pi/4}\frac{1}{\sqrt{k_0\rho}} e^{ik_0 \rho} e^{-ik_0 x^\prime \cos\alpha -ik_0 y^\prime \sin \alpha}  \nonumber \\
  &= \frac{i}{4} \sqrt{\frac{2}{\pi}}e^{-i\pi/4}\frac{1}{\sqrt{k_0\rho}} e^{ik_0\rho} e^{-ik_0 \x^\prime \cdot \boldsymbol\alpha} 
\end{align}
where we have used the fact that the Hankel function $H_0^{(1)}$ is asymptotically given by
\begin{equation}
  H_0^{(1)}(r) = \sqrt{\frac{2}{\pi r}} \exp\left( i (r - \frac{\pi}{4})\right), \qquad r \in \mathbb{R}, \quad r \gg 1.
\end{equation}
This leads to the following asymptotic form of the 2D scattered wave solution
\begin{equation}\label{eq:farfield}
\begin{split}
  u(\rho,\alpha) 
  &= \frac{i}{4}\sqrt{\frac{2}{\pi}}e^{-i\pi/4} \frac{e^{ik_0\rho}}{\sqrt{k_0\rho}} \int_\Omega e^{-ik_0 \x^\prime \cdot \boldsymbol \alpha} g(\x^\prime) \, d\x^\prime,
\end{split}
\end{equation}
for $\rho \to \infty$. The above expression is called the 2D far field wave pattern of $u$, with the integral being denoted as the far field (amplitude) map
\begin{equation} \label{eq:FFmap}
	F(\boldsymbol \alpha) = \int_\Omega e^{-ik_0 \x^\prime \cdot \boldsymbol \alpha} g(\x^\prime) \, d\x^\prime.
\end{equation}
The value of the integral depends only on the direction
$\boldsymbol\alpha$ (or, in 2D, on the angle $\alpha$) and the 
wavenumber $k_0$. Expression (\ref{eq:farfield}) readily extends to the
$d$-dimensional case, where it holds more generally that
\begin{equation} \label{eq:FFusc}
	\lim_{\rho \to \infty} u(\rho,\boldsymbol\alpha) = D(\rho) F(\boldsymbol\alpha), \qquad \boldsymbol\alpha \in \mathbb{R}^{d},
\end{equation}
for a function $D(\rho)$ which is known explicitly and a far field map $F(\boldsymbol\alpha)$ given by (\ref{eq:FFmap}).
Note that this far field map is in fact a Fourier integral of the function $g$. 

\subsection{Comments on the classical derivation} \label{subsec:comments}

It is clear from the above derivation that the calculation of the far field wave
pattern of a scattered wave consists of two main steps.  First,
one has to solve Helmholtz equation (\ref{eq:helmholtz_in}) 
with a spatially dependent wavenumber on a numerical box with 
absorbing boundary conditions. Once the numerical solution is obtained, it
is followed by the calculation of a Fourier integral (\ref{eq:FFmap})
over the aforementioned numerical domain. The main computational
bottleneck of the calculation generally lies within the first step,
since this requires an efficient and computationally inexpensive 
method for the solution of the indefinite Helmholtz 
system with absorbing boundary conditions.

The statement of the far field map presented in this section relies on the fact that the object of interest, represented by the function $\chi$, is compactly supported. In particular, this is used when computing the numerical solution $u^N$ to equation (\ref{eq:helmholtz_in}) on a bounded numerical box that covers the support of $\chi$. The above reasoning can however be readily extended to the more general class of analytical object functions $\chi$ that vanish at infinity, i.e.\ $\chi\in V$ where $V = \{f: \mathbb{R}^d \to \mathbb{R} ~\textrm{analytical}~ | ~ \forall \varepsilon > 0, \, \exists K\subset \mathbb{R}^d ~\textrm{compact}, \, \forall x \in \mathbb{R}^d \setminus K : | f(x) | < \varepsilon\}$. Indeed, due to the existence of smooth bump functions \cite{johnson2007saddle,mead1973convergence}, functions with compact support can be shown to be dense within the space of functions that vanish at infinity. Consequently, every analytical function $\chi \in V$ can be arbitrarily closely approximated by a series of compactly supported functions $\{\chi_n\}_n$. This in turn implies that the corresponding solutions $\{u^N_n\}_n$ on a limited computational box can be arbitrarily close to the solution of the Helmholtz equation generated with the analytical object of interest $\chi \in V$. Intuitively, this means that if $\chi$ is analytical but sufficiently small everywhere outside $O$, the computational domain may be retricted to a numerical box covering $O$ as if $\chi$ was compactly supported. Hence, the far field map (\ref{eq:FFmap}) is well-defined for analytical functions $\chi$ that vanish at infinity. This observation will prove particularly useful in the next section.

%
%
%
%

\section{Far field integral: complex contour formulation}\label{sec:Far}
In this section we illustrate how the far field integral (\ref{eq:FFmap}) can be reformulated as an integral over a complex contour. This new insight is consequently shown to be particularly useful with respect to the numerical computation, as it allows replacing equation (\ref{eq:helmholtz_in}) by a damped Helmholtz equation.

\begin{figure}
\begin{center}
\includegraphics[width=0.7\textwidth]{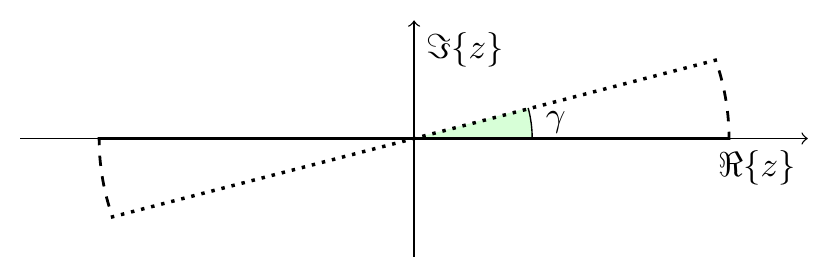} \vspace{-0.3cm}
\caption{Schematic representation of the complex contour for the far field integral calculation illustrated in 1D. The full line represents the real domain $\Omega$, the dotted and dashed lines represent the subareas $Z_1=\{ \x e^{i\gamma} : x \in \Omega \subset \mathbb{R} \}$ and $Z_2= \{b e^{i\theta} : b \in \partial\Omega,~ \theta \in [0,\gamma]\}$ of the complex contour respectively.}
\label{fig:compcont}
\end{center}
\end{figure}

\subsection{Reformulation to a complex contour} \label{subsec:reformulation}

The far field integral (\ref{eq:FFmap}) can be split into a sum of two
contributions: $F(\boldsymbol \alpha) = I_1 + I_2$, with
\begin{eqnarray} 
 	I_1 &=  k_0^2 \int_\Omega e^{-i k_0 \x \cdot \boldsymbol \alpha } \chi(\x)u_{in}(\x) d\x, \label{eq:i1} \\
	I_2 &= k_0^2 \int_\Omega e^{-i k_0 \x \cdot \boldsymbol \alpha } \chi(\x) u^N(\x) d\x.\label{eq:i2}
\end{eqnarray}
Calculation of the integral $I_1$ is generally easy, since it only
requires the expression for the incoming wave which is known
analytically. The second integral however requires the solution of the
Helmholtz equation on the domain $\Omega$ (numerical box), which is notoriously 
hard to obtain using iterative methods.
Hence, the calculation of the scattered wave $u^N$ forms a major bottleneck 
in the traditional calculation of the far field map.

However, if both $u$ and $\chi$ are analytical functions the integral can be calculated over a complex
contour rather than the real domain as follows. Let us define a complex contour along the rotated real domain $Z_1 = \{\z \in \mathbb{C} \,|\, \z = \x e^{i\gamma} : \x \in \Omega \}$, where $\gamma$ is a fixed rotation angle, followed by the curved segment $Z_2 = \{\z \in \mathbb{C} \,|\, \z = \bb e^{i\theta} : \bb \in \partial\Omega,~ 0 \leq \theta \leq \gamma\}$, as presented schematically on Figure \ref{fig:compcont} for a 1D domain. The extension of the domain to multiple dimensions is straightforward, see \cite{bayliss1985accuracy,lairdpreconditioned,reps2010indefinite}. Integral $I_2$ can then be written as
\begin{equation}
  I_2 = k_0^2 \int_{Z_1} e^{-i k_0 \z \cdot \boldsymbol \alpha }\chi(\z)u^N(\z) d\z + k_0^2 \int_{Z_2} e^{-i k_0 \z \cdot \boldsymbol \alpha }\chi(\z)u^N(\z)  d\z.
\end{equation}
The second term in the above expression vanishes, since the function $\chi$ is per definition zero everywhere outside the object of interest $O$, thus notably in all points $\z \in Z_2$. Hence, we obtain
\begin{equation} \label{eq:I2complex}
  I_2 = k_0^2 \int_{Z_1} e^{-i k_0 \z \cdot \boldsymbol \alpha }\chi(\z)u^N(\z) d\z = k_0^2 \int_{\Omega} e^{-i k_0 e^{i\gamma} \x \cdot \boldsymbol \alpha }\chi(\x e^{i\gamma})u^N(\x e^{i\gamma}) e^{i\gamma} \, d\x.
\end{equation} 
Note that for $0 < \gamma < \pi/2$ the exponential of $\x e^{i\gamma}$ is increasing in all
directions. At the same time the scattered wave solution $u^N$, which consists of outgoing
waves on the complex domain $Z_1$, is decaying in all directions. Additionally, the function $\chi$ is presumed to have a
bounded support (or vanish at infinity, see Section \ref{subsec:comments}), making the above integral computable on a limited numerical domain.

Formulation (\ref{eq:I2complex}) of the integral $I_2$ is theoretically equivalent to the 
original formulation (\ref{eq:i2}), since both formulations
result in the same value for the integral. However, the reformulation 
to the complex contour provides a significant advantage from a computational point of 
view. Indeed, (\ref{eq:I2complex}) indicates that
the far field map can (at least partially) be computed over the full
complex contour $Z_1$, that is, a rotation of the original real-valued domain
$\Omega$ over an angle $\gamma$ in all spatial dimensions.
Contrary to the original integral formulation (\ref{eq:i2}), which required the scattered wave $u^N$ 
evaluated along the real domain, formulation (\ref{eq:I2complex}) now requires the 
scattered wave $u^N$ along the complex contour. Consequently, for our new approach, the 
first step in the far field map calculation consists of solving 
the Helmholtz equation (\ref{eq:helmholtz_in}) on a complex contour, i.e.
\begin{equation} \label{eq:helmholtz_cc}
  \left(-\Delta - k^2(\z)\right) u(\z) = k_0^2 \chi(\z)u_{in}(\z), \qquad \z \in Z_1,
\end{equation}
which is known to be much easier to solve than the original real-valued Helmholtz equation.
Classical multigrid methods allow for a fast solution of damped Helmholtz systems, see \cite{erlangga2006novel}.

The above integral reformulation (\ref{eq:I2complex}) implies no restrictions on the 
value of the rotation angle. Indeed, the rotation angle $\gamma$
can in principle be chosen anywhere in the interval $[0, \pi/2]$, as the 
equivalence between formulation (\ref{eq:i2}) and (\ref{eq:I2complex}) theoretically holds 
for any angle $\gamma$. However, it will be shown in Section 
\ref{subsec:rot_angle} that in practice, a lower bound on the rotation angle $\gamma$ is implied by 
the numerical computation of the scattered wave $u^N$ in the first step of the 
far field map computation. Note that, on the other hand, the principal incentive for not choosing $\gamma$ excessively large 
is the stability of the numerical integration scheme. When the angle $\gamma$ is very large, 
the scattered wave is heavily damped in all regions remote from the origin, and the integration scheme 
typically requires additional evaluation points near the origin, cf.\ \cite{huybrechs2006evaluation}. 

%
%

\begin{figure}
\begin{center}
\includegraphics[width=0.7\textwidth]{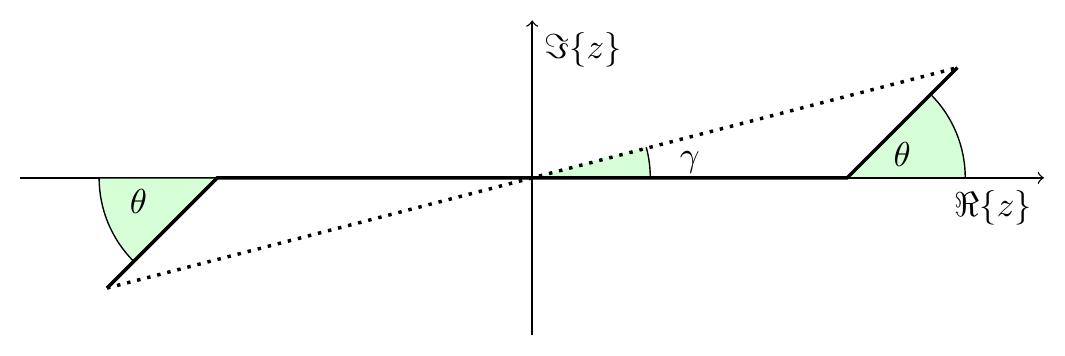} \vspace{-0.3cm}
\caption{Schematic representation of a real grid with absorbing boundaries with ECS angle $\theta$ (solid line) vs.\ a full complex grid with rotation angle $\gamma$ (dotted line), illustrated in 1D. The solid line is the classical domain for the real-valued solution $u^N(\x)$ to (\ref{eq:helmholtz_in}), whereas the dotted line represents the domain $Z_1$ for the complex-valued solution $u^N(\z)$. \label{fig:contour}}
\label{fig:doubleecs}
\end{center}
\end{figure}

\subsection{Solving the Helmholtz equation on a complex domain} \label{subsec:solving}

We now show that the formulation of the Helmholtz problem 
on a complex rotated domain like $Z_1$ is very similar to 
a complex shifted Laplacian system \cite{erlangga2006novel}. 
Indeed, both formulations can be shown to be generally 
equivalent, and hence are equally efficiently solvable. 
Consider a Helmholtz problem on a complex rotated 
grid of the form (\ref{eq:helmholtz_cc})
\begin{equation} \label{eq:HHCC}
  \left(-\Delta - k^2(\z)\right)u(\z) = b(\z), \qquad \z \in Z_1.
\end{equation}
This equation is discretized using finite differences on a 
$d$-dimensional Cartesian grid with a complex valued grid 
distance $\tilde{h} = e^{i \gamma} h$ (with $h \in \mathbb{R}$) 
in every spatial dimension, yielding a linear system of the form
\begin{equation}\label{eq:csg}
  -\left(\frac{1}{e^{i \varphi} h^2} L  + k^2\right) u_h = b_h,
\end{equation}
with $\varphi = 2\gamma$ and where $L$ is the matrix operator expressing 
the stencil structure of the Laplacian. For example, discretization of the 
2D Laplacian using second order finite differences yields 
$L = \text{kron}(I, \text{diag}(-1,2,-1) + \text{kron}(\text{diag}(-1,2,-1),I)$, 
where the size of $L$ intrinsically depends on $h$. Rewriting the rotation as 
$\exp{(i \varphi)} = (\alpha + i\beta)$, and multiplying both sides in 
\eqref{eq:csg} by $(\alpha+i\beta)$, we obtain the equivalent system
\begin{equation}\label{eq:csl}
  -\left(\frac{1}{h^2} L  + (\alpha+i\beta)k^2\right) u_h = (\alpha+\beta i) \, b_h.
\end{equation}
The left-hand side matrix operator of this equation is a discretization of 
the complex shifted Helmholtz operator $-\Delta-(\alpha+\beta i)k^2$. 
We momentarily assume $k$ to be real-valued; an analogous 
argument can be used for complex-valued wavenumbers 
under proper conditions.  Hence, (\ref{eq:csl}) is the discrete 
representation of a CSL system, which is known to be solvable using 
multigrid given a sufficiently large complex shift $\beta$ 
\cite{erlangga2008multilevel}. Note that in the CSL literature,
the real part of the shift factor is commonly chosen as $\alpha = 1$, 
yielding a shift of the form $-(1+i\beta) k^2$ \cite{erlangga2006novel}.
We expound on the relation between the complex shift $\beta$ and 
the rotation angle $\gamma$ in Section \ref{subsec:rot_angle}.

The scheme given by (\ref{eq:csg}) is known as Complex Stretched Grid (CSG), 
and was shown in \cite{reps2010indefinite} to hold very similar properties 
with respect to multigrid convergence compared to the complex shifted 
Laplacian system. Moreover, it is shown above that the CSL and CSG schemes 
are generally equivalent, and thus can be solved equally efficiently using 
a multigrid method.

\subsection{On the rotation angle}\label{subsec:rot_angle}

The choice of a sufficiently large complex shift parameter $\beta$ is
vital to the stability of the multigrid solution method for problem \eqref{eq:csl}, or equivalently \eqref{eq:csg}. 
When the shift parameter is chosen below a certain minimal value denoted as $\beta_{min}$, the multigrid scheme tends to be unstable and convergence will break down.
A typical rule of thumb for the choice of the complex shift suggested in the CSL
literature is $\beta \geq 0.5$ \cite{cools2012local,erlangga2006novel}. 
Note that the lower limit $\beta_{min} = 0.5$ is based on a multigrid V-cycle with standard weighted
Jacobi or Gauss-Seidel smoothing. However, more advanced iterative techniques 
like ILU \cite{umetani2009multigrid} or GMRES \cite{calandra2011two} may alternatively be used as 
a smoother substitute in the multigrid solver.

The requirement of a minimal shift for multigrid stability on the CSL problem can be directly
translated into a minimal rotation angle $\gamma$ for the complex scaled system. 
Writing the complex shift in polar notation
\begin{equation}
(1+i\beta) = \rho \exp(i\varphi),
\end{equation}
where $\rho=\sqrt{1+\beta^2}$ and $\varphi=\arctan\beta$, one readily obtains
\begin{equation}
  \tilde{h} = \sqrt{1+i\beta} \, h = \sqrt{\rho} \exp(i\varphi/2) \, h.
\end{equation}
Denoting the minimal value of the complex shift by $\beta_{min}$, 
the rotation angle $\gamma = \varphi/2$ must satisfy
\begin{equation} \label{eq:mingammajac}
   \gamma \geq \frac{\arctan(\beta_{min})}{2}.
\end{equation}
For the rule of thumb $\beta_{min} = 0.5$ stated above, this implies that $\gamma \geq 0.23 \approx 13^\circ$.
Note that when substituting the standard stationary multigrid relaxation scheme ($\omega$-Jacobi, Gauss-Seidel) 
by a more robust iterative scheme like GMRES(m) (with $m = 2$ or $3$), the rotation angle $\gamma$ is typically chosen at around $10^\circ$, resulting in a stable multigrid scheme.

In this paper we have chosen to link the grid rotation angle $\gamma$
to the standard ECS absorbing layer angle $\theta$, as shown in Figure
\ref{fig:doubleecs}. This is in no way imperious for the functionality
of the method, but it appears quite naturally from the fact that both
angles perturb (part of) the grid into the complex plane. Supposing the
ECS boundary layer measures one quarter of the length of the entire
real domain in every spatial dimension, which is a common choice, we
readily derive the following relation between the rotation angle $\gamma$
and the ECS angle $\theta$:
\begin{equation} \label{eq:thetagamma}
	\gamma = \arctan\left(\frac{\sin\theta}{2 + \cos\theta}\right).
\end{equation}
Table \ref{tab:gammatheta} shows some standard values of the ECS angle
$\theta$ and corresponding $\gamma$ values according to (\ref{eq:thetagamma}). 
Following (\ref{eq:mingammajac}), we note that for a multigrid scheme with
$\omega$-Jacobi or Gauss-Seidel smoothing to be stable, $\theta$
should be chosen no smaller than $\pi/4$. Using the more efficient GMRES(3) 
method as a smoother replacement, an ECS angle around $\theta = \pi/6$ 
suffices to guarantee stable multigrid convergence.

On the other hand, taking $\gamma$ too large has some drawbacks which are
reflected in the accuracy of the numerical scheme for the calculation of 
the far-field integral. For $\gamma > 0$, the solution $u^N(\z)$ is an 
exponentially decaying function towards the boundary. At the same time,
the function $e^{-i k_0  \z \cdot \mathbf{\alpha}}$ is exponentially 
growing towards the boundary. Their product, which appears in the integral, 
remains bounded. However, when $\gamma$ is chosen to be very large, the 
difference in magnitude between the two integrand factors can affect the 
accuracy of the numerical integral. Indeed, when multiplying two floating 
point numbers, one being extremely small and the other extremely large, 
and each is represented up to machine precision, the numerical error on this 
multiplication can be very large. These considerations somewhat limits the choice 
of the rotation angle for the complex contour.

\begin{table}[t]
\centering
\begin{tabular}{| c || c | c | c | c | c | c |}
\hline
	$\theta$ (rad.) & $\pi/8$ & $\pi/7$ & $\pi/6$ & $\pi/5$ & $\pi/4$ & $\pi/3$ \\
\hline
	$\gamma$ (deg.) & ~$7.5^\circ$ & ~$8.5^\circ$ & ~$9.9^\circ$ & ~$11.8^\circ$ & ~$14.6^\circ$ & ~$19.1^\circ$ \\
\hline
\end{tabular}
\vspace{0.2cm}
\caption{ECS angle $\theta$ and corresponding rotation angle $\gamma$ for the full complex grid. Values based on (\ref{eq:thetagamma}).}
\label{tab:gammatheta}
\end{table}

%
%
%
%

\section{Numerical results for 2D and 3D Helmholtz problems}\label{sec:Numerical}

In this section, we validate the theoretical result presented above
by a number of numerical experiments in both two and three spatial
dimensions. We make use of the known efficiency of multigrid in solving damped Helmholtz 
equations to compute the solutions to the Helmholtz systems required in the complex valued far field integral.

\begin{figure}[t!] 
\begin{center}
\includegraphics[width=0.48\textwidth]{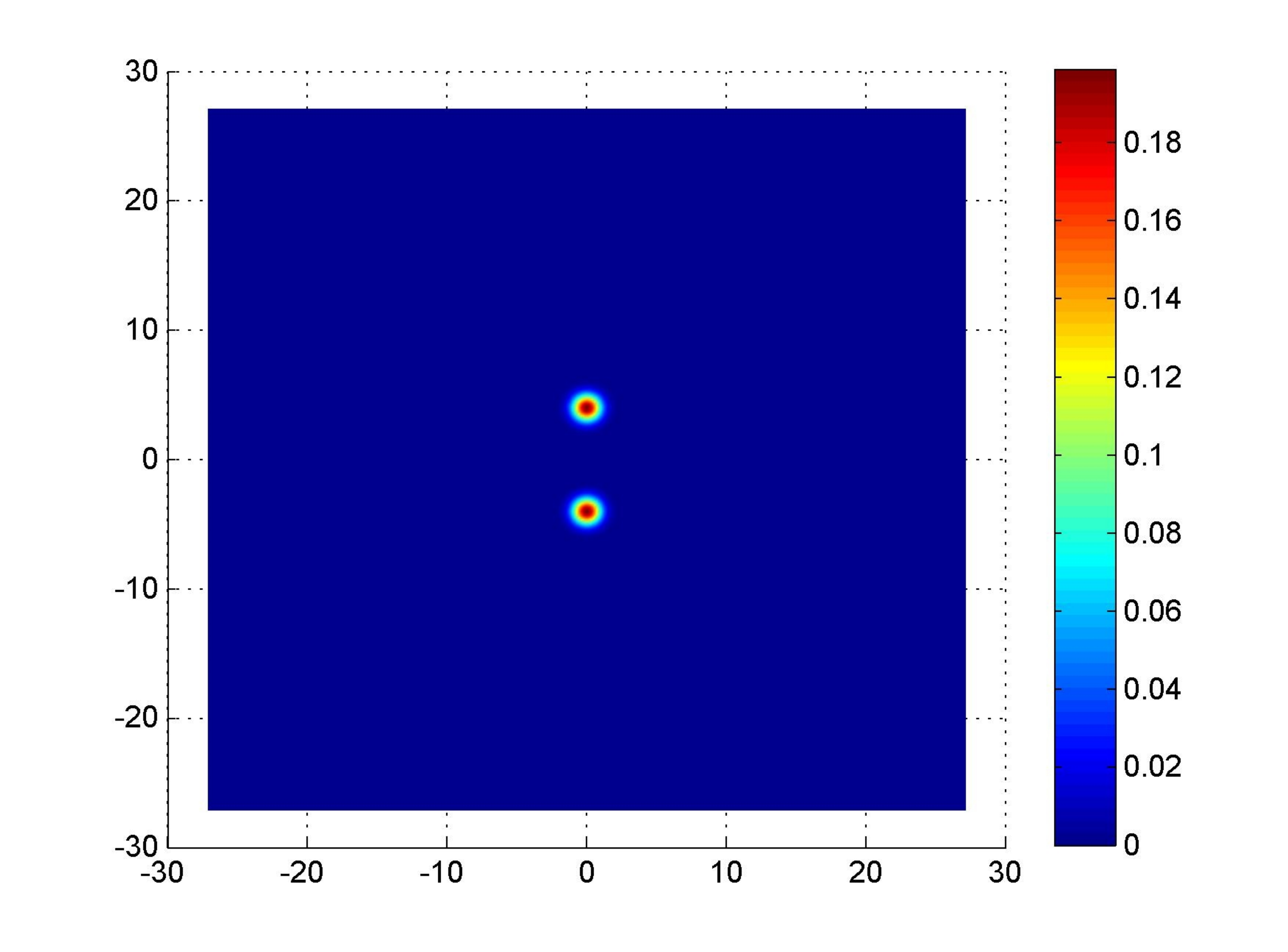}
\includegraphics[width=0.48\textwidth]{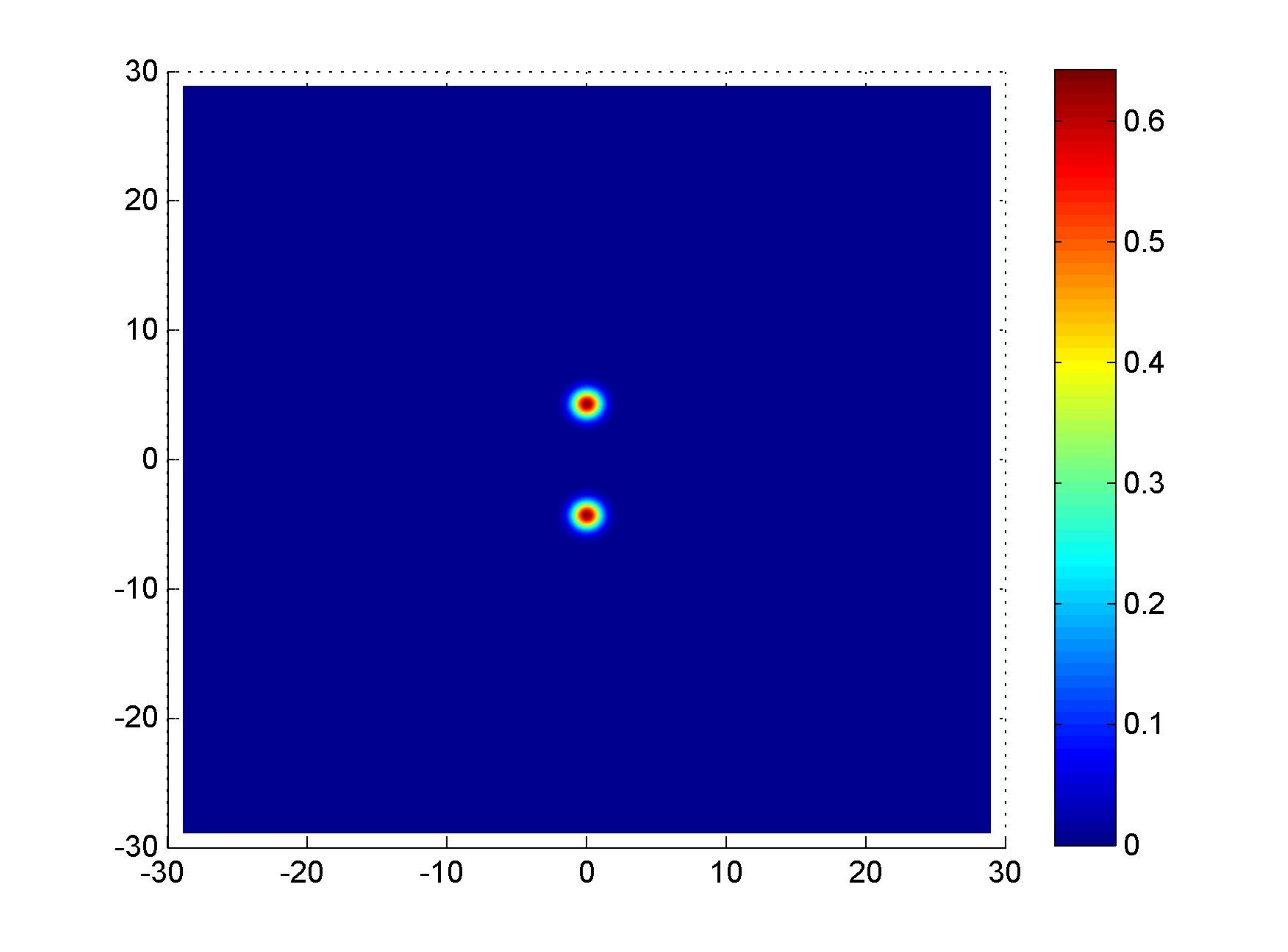}\\ 		\vspace{0.0cm}
\includegraphics[width=0.48\textwidth]{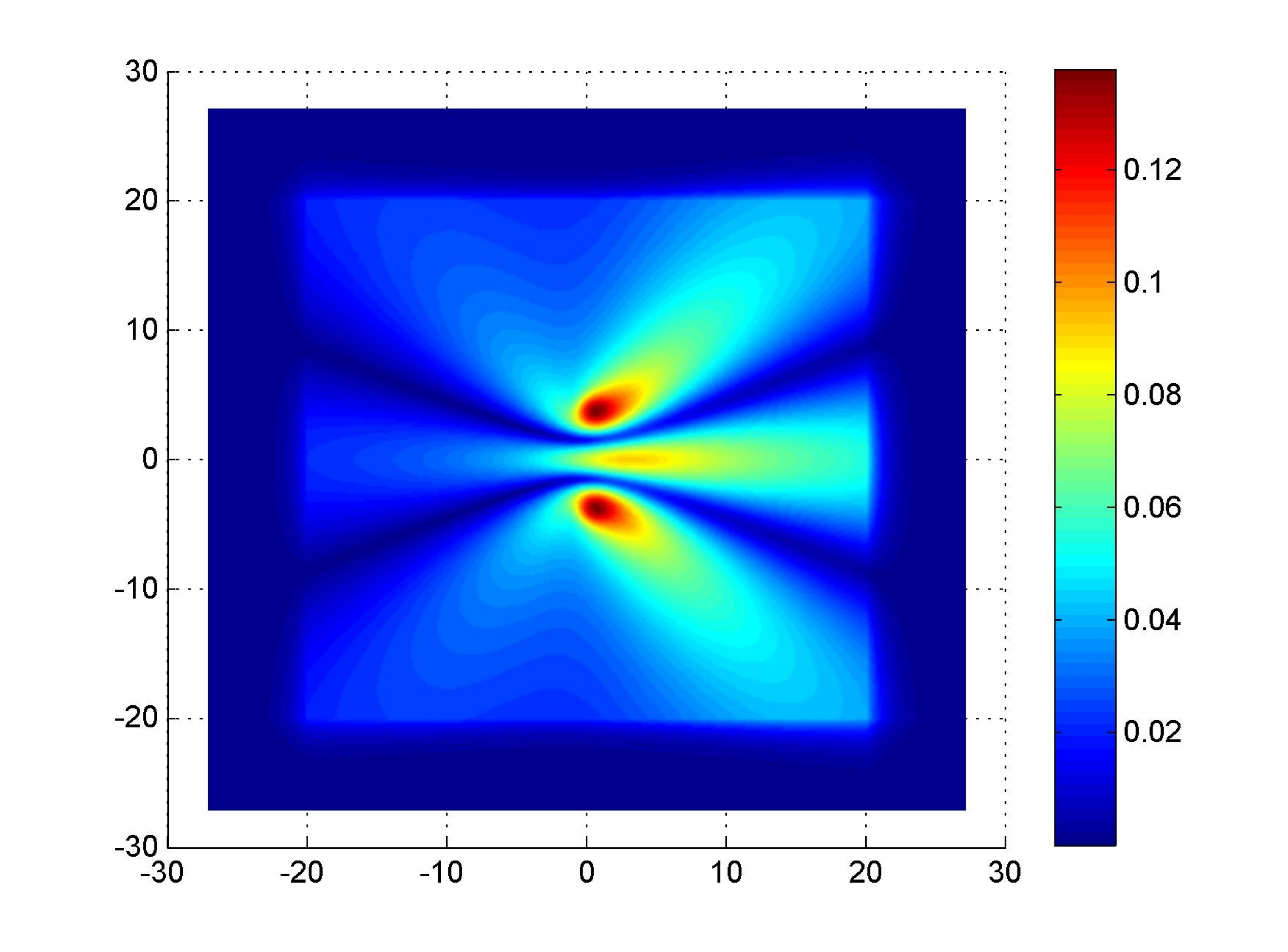}
\includegraphics[width=0.48\textwidth]{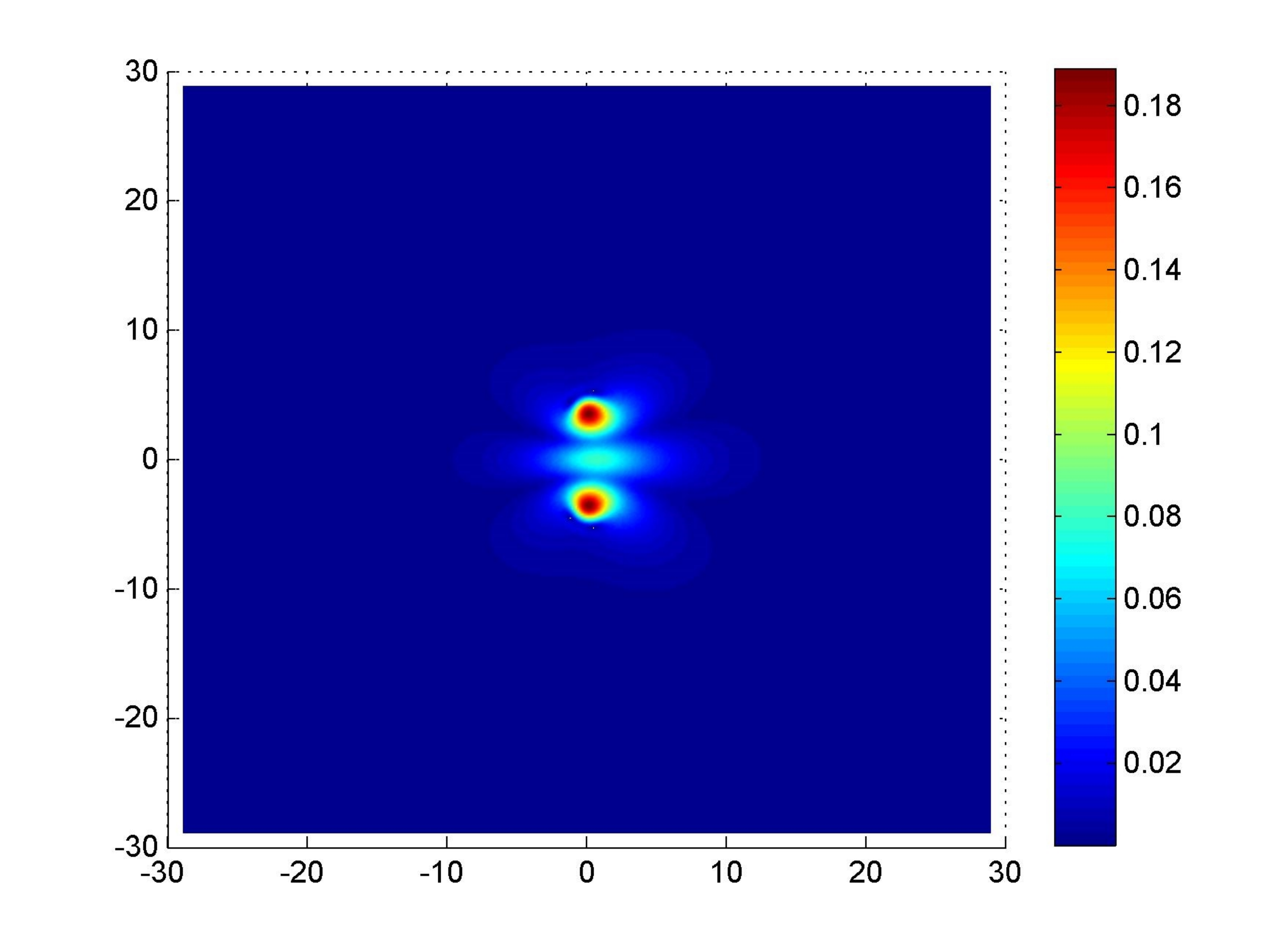}\\ 		\vspace{0.0cm}
\includegraphics[width=0.48\textwidth]{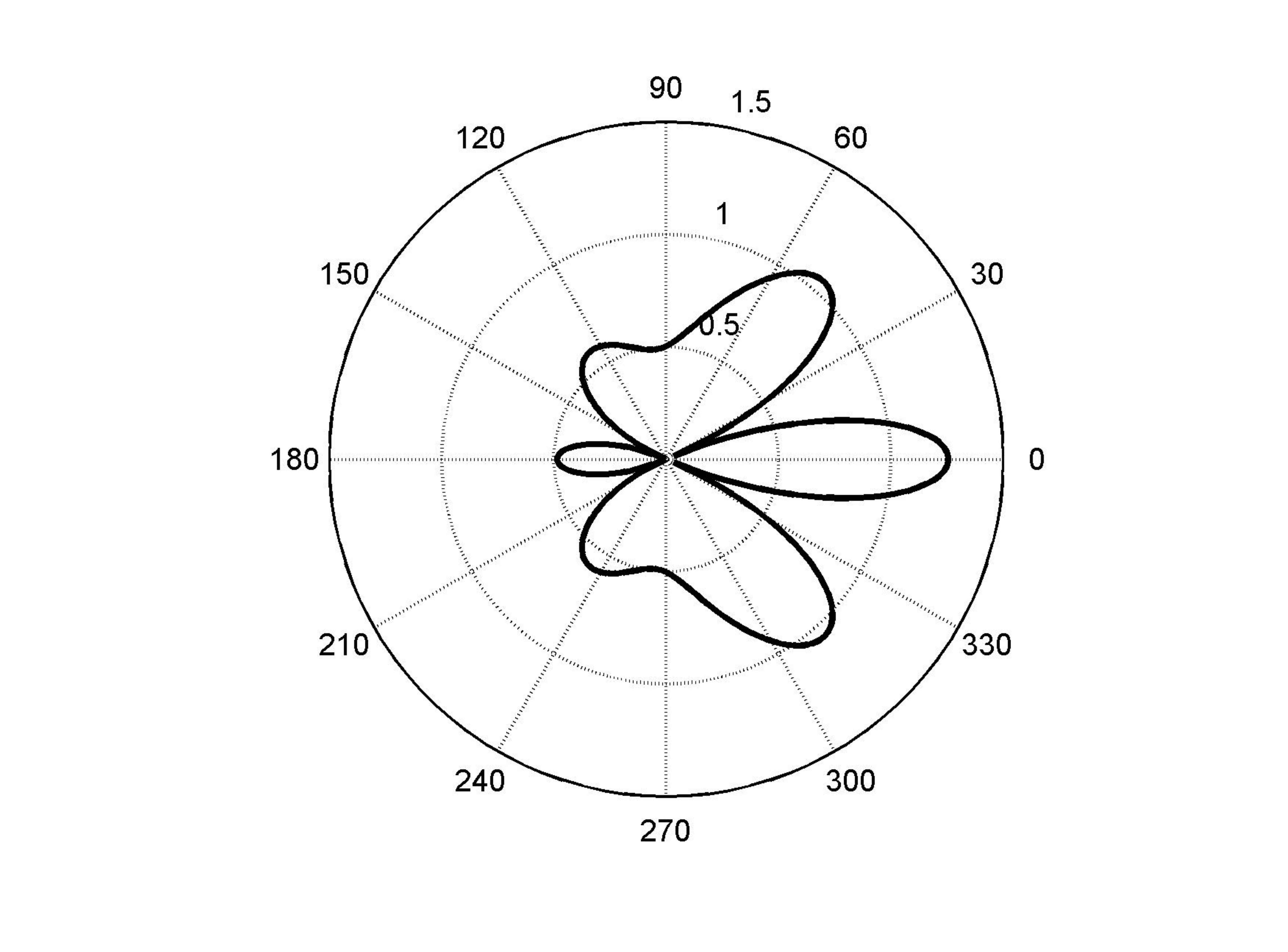}
\includegraphics[width=0.48\textwidth]{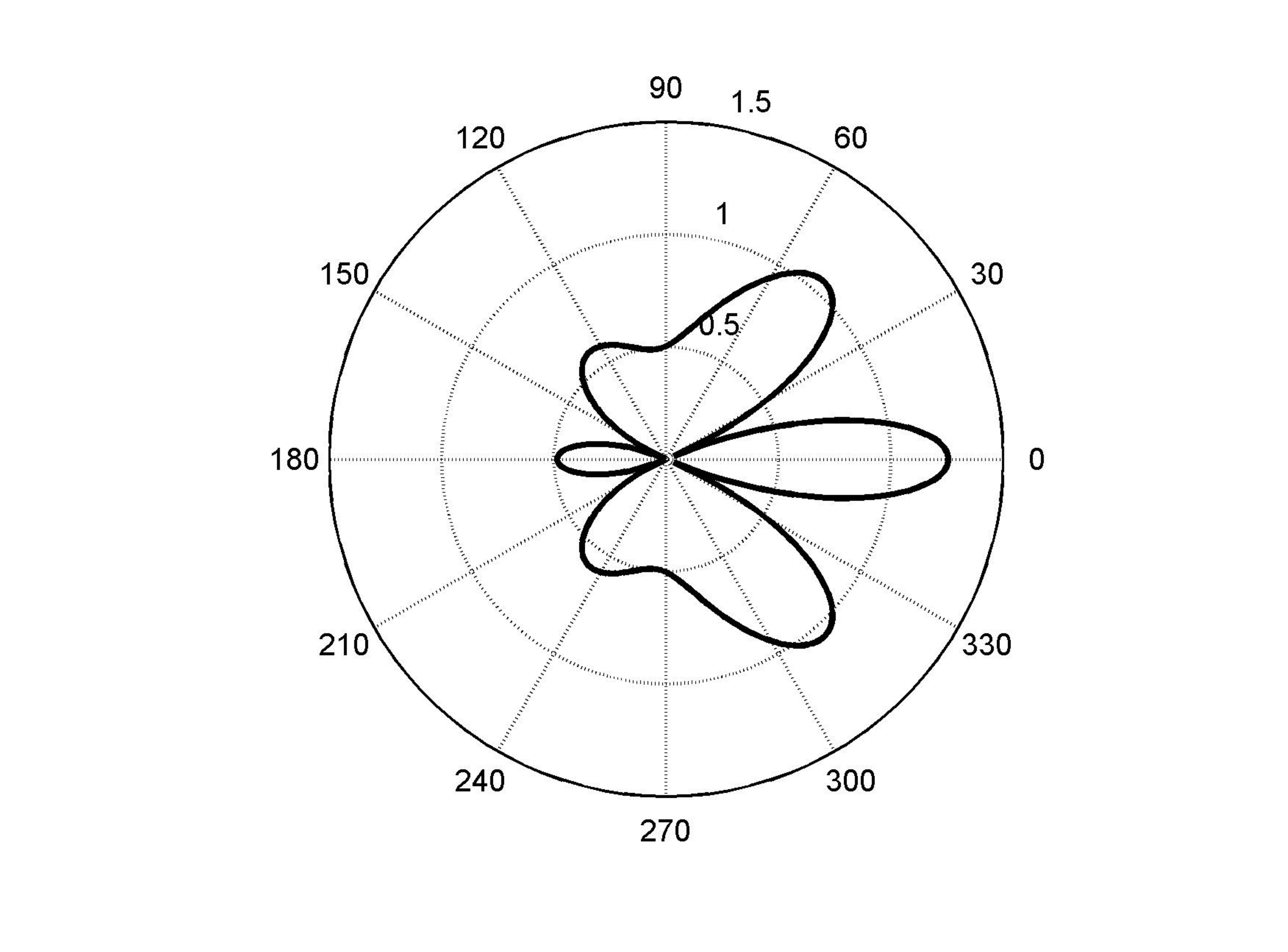} 	  \vspace{-0.4cm}
\caption{Comparison between the classical real-valued far field map calculation method (left) and the new complex-valued approach (right). Top: 2D object of interest $|\chi|$ given by (\ref{eq:chi2D}). Mid: solutions to the Helmholtz problem (\ref{eq:helmholtz_in}) (in modulus) on a $n_x \times n_y = 256 \times 256$ grid. Left: real-valued problem with double ECS contour with $\theta = \pi/4$ solved using LU factorization. Right: complex-valued problem solved using a series of multigrid V-cycles with $\omega$-Jacobi smoother on the corresponding full complex contour with $\gamma = 14.6^\circ$ up to a residual reduction tolerance of \texttt{1e-6}. Bottom: resulting 2D Far field map $F(\alpha)$ for both approaches, calculated following (\ref{eq:FFmap}). Normalized errors with respect to a $n_x \times n_y = 1024\times1024$ grid benchmark Far field map solution $F_{ex}(\alpha)$: $\|F_{re}-F_{ex}\|_2/\|F_{ex}\|_2 =$ \texttt{9.37e-5} (left), $\|F_{co}-F_{ex}\|_2/\|F_{ex}\|_2 =$ \texttt{1.39e-4} (right), $\|F_{re}-F_{co}\|_2/\|F_{ex}\|_2 =$ \texttt{1.77e-4}.}
\label{fig:results2D}
\end{center}
\end{figure}

\begin{figure}[t!] 
\begin{center}
\includegraphics[width=0.48\textwidth]{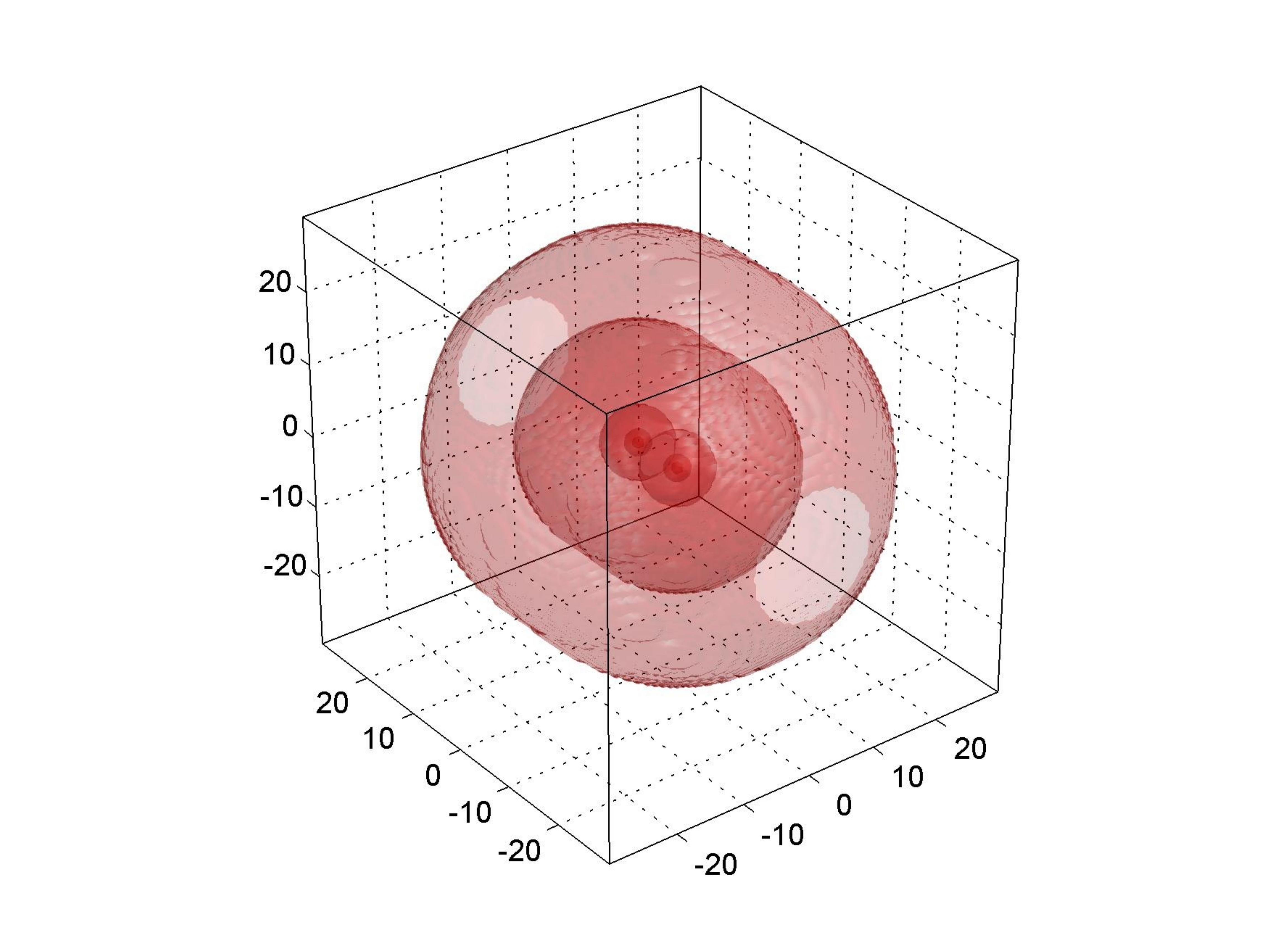}
\includegraphics[width=0.48\textwidth]{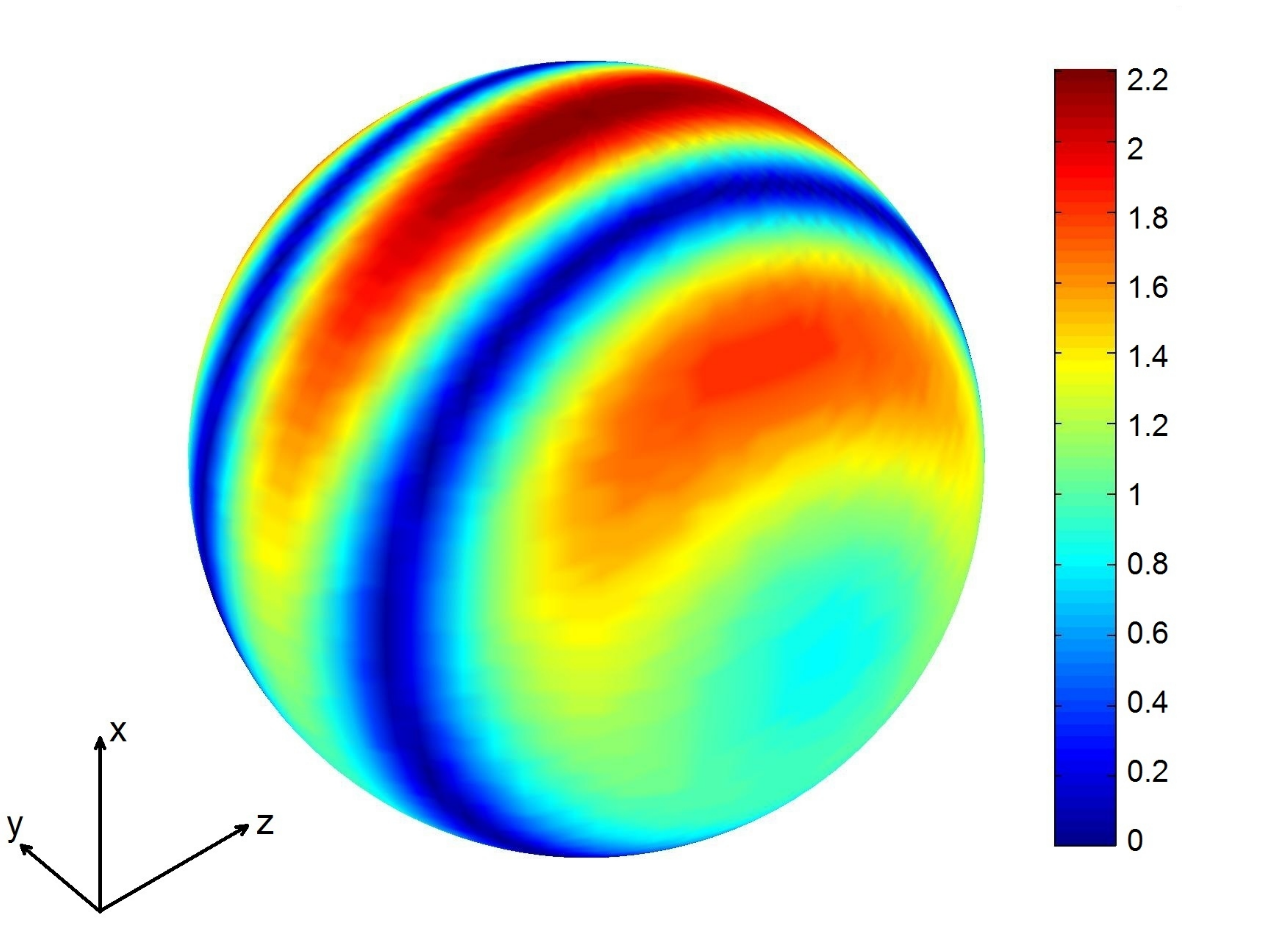}
\caption{Left: 3D object of interest $|\chi(\x)|$ given by (\ref{eq:chi3D}). Shown are the $|\chi(\x)| = c$ isosurfaces for $c =$ \texttt{1e-1}, \texttt{1e-2}, \texttt{1e-10}, \texttt{1e-100} and \texttt{1e-300}. Right: 3D Far field map, resulting from Helmholtz problem (\ref{eq:helmholtz_in}) with $k_0 = 1$ solved on a $n_x \times n_y \times n_z = 64 \times 64 \times 64$ full complex grid with $\theta = \pi/6$ ($\gamma \approx 9.9^\circ$) using a series of multigrid V-cycles with GMRES(3) smoother up to residual reduction tolerance \texttt{1e-6}.}
\label{fig:ffmap3D}
\end{center}
\end{figure}

The model problem used throughout this section is a Helmholtz equation
of the form (\ref{eq:helmholtz_in}) with $k^2(\x) = k_0^2 (1 +
\chi(\x))$. The equation is discretized on a $n^d$-point uniform Cartesian mesh
covering a square numerical domain $\Omega = [-20,20]^d$ using second
order finite differences. Note that the use of a different discretization 
scheme would not fundamentally affect the results presented in this work.
In the 2D case the space-dependent wavenumber is defined as
\begin{equation} \label{eq:chi2D}
	k_0^2 \, \chi(x,y) = - 1/5 \left( e^{-(x^2+(y-4)^2)} + e^{-(x^2+(y+4)^2)} \right), \quad (x,y) \in [-20,20]^2,
\end{equation}
i.e.\ the object of interest takes the form of two circular point-like
objects with mass concentrated at the Cartesian coordinates $(0,-4)$
and $(0,4)$, see Figure \ref{fig:results2D} (top panel). For the 3D model
problem, the following straightforward extension of the object is used
\begin{equation} \label{eq:chi3D}
	k_0^2 \, \chi(x,y,z) = - 1/5 \left( e^{-(x^2+(y-4)^2+z^2)} + e^{-(x^2+(y+4)^2+z^2)} \right), \quad (x,y,z) \in [-20,20]^3,
\end{equation}
representing two spherical point-like objects in 3D space, see Figure
\ref{fig:ffmap3D} (left panel). The incoming wave scattering at the given object
is defined by
\begin{equation}
	u_{in}(\x) = e^{i k_0 \boldsymbol\eta \cdot \x}, \quad \x \in \Omega,
\end{equation}
where $\boldsymbol\eta$ is the unit vector in the $x$-direction.

\begin{table}[t]
\vspace{0.5cm}
\centering
\begin{tabular}{| c | c | c c c c c |}
\hline
 & \multirow{2}{*}{{\footnotesize$n_x \times n_y \times n_z$}} & \multirow{2}{*}{~~~~$16^3$~~~~} & \multirow{2}{*}{~~~~$32^3$~~~~} & \multirow{2}{*}{~~~~$64^3$~~~~} & \multirow{2}{*}{~~~~$128^3$~~~~} & \multirow{2}{*}{~~~~$256^3$~~~~} \\
 &  &  &  &  &  & \\
\hline
\multirow{10}{*}{~$k_0$~} & \multirow{2}{*}{1/4} 	& \textbf{10} (10) 		& \textbf{9} (59) 		& \textbf{9} (560) 		& \textbf{9} (4456) 	& \textbf{9} (35165)		\\
 													&    										& 0.24								& 0.20		 						& 0.21								& 0.20 								& 0.20									\\
 													&	\multirow{2}{*}{1/2} 	& 12 (12)							& \textbf{10} (63) 		& \textbf{10} (611) 	& \textbf{10} (4937)	& \textbf{9} (35405)		\\
 													&    										& 0.31								& 0.24		 						& 0.22								& 0.23 								& 0.21									\\
 													&	\multirow{2}{*}{1}   	& 7  (8) 							& 13 (83) 						& \textbf{11} (691)		& \textbf{10} (4899)	& \textbf{10} (38975)		\\
 													&    										& 0.13								& 0.32		 						& 0.27								& 0.24 								& 0.24									\\
 													&	\multirow{2}{*}{2}   	& 2  (4) 							& 8 (54) 							& 13 (809) 						& \textbf{11} (5418)	& \textbf{10} (38051)		\\
 													&    										& 0.01								& 0.14		 						& 0.33								& 0.27 								& 0.24									\\
 													&	\multirow{2}{*}{4}   	& 1  (3) 							& 2 (17) 							& 7 (457) 						& 13 (6337)						& \textbf{11} (41848)		\\
 													&	   										& 0.01								& 0.01		 						& 0.12								& 0.33 								& 0.26									\\	
\hline
\end{tabular}
\vspace{0.2cm}
\caption{3D Helmholtz problem (\ref{eq:helmholtz_in}) solved on a full complex grid with $\theta = \pi/6$ $(\gamma=9.9^\circ)$ using a series of multigrid V(1,1)-cycles with GMRES(3) smoother up to residual reduction tolerance \texttt{1e-6}. Displayed are the number of V-cycle iterations, number of work units and average convergence factor for various wavenumbers $k_0$ and different discretizations. 1 WU is calibrated as the cost of 1 V(1,1)-cycle on the $16^3$-points grid $k_0 = 1/4$ problem. Discretizations respecting the $k_0 h< 0.625$ criterion for a minimum of 10 grid points per wavelength are indicated by a bold typesetting.}
\label{tab:Vkh3D}
\end{table}

Figure \ref{fig:results2D} validates the equivalence between the classical 
real-valued far field map integral from Section \ref{sec:Helmholtz} and the new complex contour formulation
presented in Section \ref{sec:Far}. The 2D Helmholtz model problem with
wavenumber given by (\ref{eq:chi2D}) and $k_0 = 1$ is solved for $u^N$ using
respectively a standard LU factorization method on the real domain
$\Omega$ with ECS complex boundary layers ($\theta = \pi/4$) along the
domain boundary $\partial\Omega$, and a series of multigrid
V(1,1)-cycles with $\omega$-Jacobi smoothing on the full complex
domain ($\gamma = 14.6^{\circ} $) with a residual reduction tolerance
of \texttt{1e-6}. The standard multigrid intergrid operators used 
in this work are bilinear interpolation and full weighting
restriction. The moduli of the wavenumber function $\chi$ (top) and
the resulting solution $u^N$ (mid) are shown on Figure
\ref{fig:results2D} for both methods. Note how the solution $u^N$ on
the full complex contour (right) is indeed heavily damped compared to the
solution on the real domain (left). Consequently, using the numerical
solution $u^N$, the 2D far field map integral (\ref{eq:FFmap})
can be calculated using any numerical integration scheme over the
real or complex domain respectively. The resulting far field map
$F(\boldsymbol\alpha)$ is shown as a function of the angle $\alpha$ on Figure
\ref{fig:results2D} (bottom). One observes that the mapping is indeed
identical when calculated over the real-valued (left) and complex-valued (right) domain, conform with
the theoretical results. The normalized difference between both far field 
computations does not exceed 0.177\% (in norm), which is effectively of the same order of 
magnitude as the normalized error. However, the computational cost of the
real-domain method for calculation of the far field map is reduced
significantly by the ability to apply multigrid to the equivalent
complex scaled problem.

\begin{table}[t]
\vspace{0.5cm}
\centering
\begin{tabular}{| c | c | c c c c c |}
\hline
 & \multirow{2}{*}{{\footnotesize$n_x \times n_y \times n_z$}} & \multirow{2}{*}{~~~~$16^3$~~~~} & \multirow{2}{*}{~~~~$32^3$~~~~} & \multirow{2}{*}{~~~~$64^3$~~~~} & \multirow{2}{*}{~~~~$128^3$~~~~} & \multirow{2}{*}{~~~~$256^3$~~~~} \\
 &  &  &  &  &  & \\
\hline
\multirow{10}{*}{~$k_0$~} & \multirow{2}{*}{1/4} 	& \textbf{8} (11) 		& \textbf{6} (52) 		& \textbf{5} (384) 		& \textbf{5} (3190) 	& \textbf{5} (25241)			\\
 													&    										& 1.93e{-9}						& 1.77e{-9}						& 2.46e{-9}						& 2.68e{-9} 					& 3.00e{-9}							\\
 													&	\multirow{2}{*}{1/2} 	& 9 (12)							& \textbf{8} (68) 		& \textbf{6} (452) 		& \textbf{6} (3392)		& \textbf{6} (30215)			\\
 													&    										& 1.27e{-8}						& 1.87e{-9}		 				& 3.37e{-9}						& 1.30e{-9} 					& 1.25e{-9}							\\
 													&	\multirow{2}{*}{1}   	& 5 (8) 							& 9 (68) 							& \textbf{8} (572)		& \textbf{7} (4013)		& \textbf{6} (30747)			\\
 													&    										& 1.33e{-8}						& 1.51e{-8}		 				& 4.07e{-9}						& 1.76e{-9} 					& 3.66e{-9}							\\
 													&	\multirow{2}{*}{2}   	& 1 (5) 							& 5 (43) 							& 9 (600) 						& \textbf{8} (4456)		& \textbf{7} (36367)		\\
 													&    										& 5.99e{-13}					& 1.18e{-8}						& 1.91e{-8}						& 5.28e{-9} 					& 2.67e{-9}						\\
 													&	\multirow{2}{*}{4}   	& 1 (4) 							& 1 (18) 							& 5 (357) 						& 9 (5038)						& \textbf{8} (39038)		\\
 													&	   										& 8.90e{-20}					& 2.86e{-13}					& 5.19e{-9}						& 1.97e{-8} 					& 4.65e{-9}						\\	
\hline
\end{tabular}
\vspace{0.2cm}
\caption{3D Helmholtz problem (\ref{eq:helmholtz_in}) solved on a full complex grid with $\theta = \pi/6$ $(\gamma = 9.9^\circ)$ using an FMG cycle with GMRES(3) smoother up to residual reduction tolerance of \texttt{1e-6}. Displayed are the number of V-cycle iterations on the designated finest grid, number of work units and resulting residual norm for various wavenumbers $k_0$ and different discretizations. 1 WU is calibrated as the cost of 1 V(1,1)-cycle on the $16^3$-points grid $k_0 = 1/4$ problem. Discretizations respecting the $k_0 h< 0.625$ criterion for a minimum of 10 grid points per wavelength indicated by a bold typesetting.}
\label{tab:HHFMG3D}
\end{table}

In Table \ref{tab:Vkh3D} convergence results are shown for the
solution of the 3D scattered wave equation (\ref{eq:helmholtz_in})
using a series of multigrid V(1,1)-cycles on various grid sizes. Note
that the multigrid method scales perfectly as a function of the number
of grid points, as doubling the number of grid points in every spatial
dimension does not increase  the number of V-cycles required to reach a
fixed residual tolerance of \texttt{1e-6}. This is a standard result
from multigrid theory. Additionally and more importantly, remarkable
$k$-scalability is measured for the multigrid solution method on the
complex contour. Indeed, the multigrid convergence factor (and thus
the corresponding work unit load required to solve the problem up to a given
tolerance) is almost fully independent of the wavenumber $k_0$, as can
be observed from the table. From a physical-numerical point
of view it is only meaningful to consider discretizations satisfying the
$k_0 h < 0.625$ criterion for a minimum of 10 grid points per
wavelength, cf.\ \cite{bayliss1985accuracy}, for which the corresponding values
are designated in Table \ref{tab:Vkh3D} by a bold typesetting.

Ultimately, the computed scattered wave solution on the 3D complex domain
can be used to calculate the far field integral
(\ref{eq:FFmap}). The resulting 3D far field mapping for the model
problem with $k_0 = 1$ is shown in Figure \ref{fig:ffmap3D}. The
left hand side panel shows an isosurface visualization of the 3D
object of interest $\chi(\x)$ given by (\ref{eq:chi3D}). On the right
panel a spherical projection of the resulting 3D far field mapping is
shown. The color hue indicates the value of the far field amplitude in
each outgoing direction.

Note that the calculation of the scattered wave solution can be optimized 
even further by considering the Full Multigrid (FMG) scheme. 
This is a nested iteration of standard V-cycles, where on each level a series of 
V(1,1)-cycles is used to approximately solve the error equation and supply a corrected initial 
guess for a finer level by interpolating the corresponding coarse grid solution.

Table \ref{tab:HHFMG3D} shows convergence results for the solution of 
the 3D scattered wave equation (\ref{eq:helmholtz_in}) using an FMG scheme. 
The setting is comparable to that of Table \ref{tab:Vkh3D}. A residual 
reduction tolerance of $10^{-6}$ is imposed for each wavenumber and at every level of the FMG cycle, yielding 
a fine $n_x \times n_y \times n_z = 256^3$ grid residual of order of magnitude $10^{-9}$. Note that the number of 
V-cycles performed on each level in the FMG cycle is decaying as a function of 
the growing grid size due to the increasingly accurate initial guess, resulting 
in a relatively small number of V-cycles (five to nine) to be performed on the finest level. 
Consequently, the number of work units (and thus the computational time) required to 
reach the designated residual reduction tolerance is significantly lower than the 
work unit load of the pure V-cycle scheme displayed in Table \ref{tab:Vkh3D}.

\begin{table}[t]
\centering
\begin{tabular}{| c | c c c c c |}
\hline
 \multirow{2}{*}{$n_x \times n_y \times n_z$} & \multirow{2}{*}{$16^3$} & \multirow{2}{*}{$32^3$} & \multirow{2}{*}{$64^3$} & \multirow{2}{*}{$128^3$} & \multirow{2}{*}{$256^3$} \\
 &  &  &  &  & \\
\hline
 CPU time			& 0.20 s.	& 0.78	s.	& 6.24 s.	& 53.3 s.	& 462	s.	\\
 $\|r\|_2$		& 3.3e-5	& 7.9e-5	& 2.7e-5	& 1.1e-5	& 4.6e-6	\\
\hline
\end{tabular}
\vspace{0.2cm}
\caption{3D Helmholtz problem (\ref{eq:helmholtz_in}) with wavenumber $k_0 = 1$ solved on a full complex grid with $\theta = \pi/6$ $(\gamma = 9.9^\circ)$ using one FMG-cycle with GMRES(3) smoother. Displayed are the CPU time (in s.) and the resulting residual norm for various discretizations. System specifications: Intel$^{{\scriptsize\textregistered}}\hspace{-0.05cm}$ Core$^{\text{\tiny TM}}$ i7-2720QM 2.20GHz CPU, 6MB Cache, 8GB RAM.}
\label{tab:Fkh3D}
\end{table}

Timing and residual results from a standard FMG sweep performing only one V(1,1)-cycle on each level on the 3D Helmholtz scattering problem with a moderate
wavenumber $k_0 = 1$ are shown in Table \ref{tab:Fkh3D} for different discretizations. Note that timings were generated using a basic non-parallelized Matlab code, using only a single thread on a simple midrange personal computer (system specifications: see caption Table \ref{tab:Fkh3D}) and taking less than 8 minutes to solve a 3D Helmholtz problem with $256$ grid points in every spatial dimension.

%
%
%
%

\section{Application to Schr\"odinger equations} \label{sec:Schrodinger}

This section illustrates the application of the proposed complex contour method to 
Schr\"odinger equations that are used to describe quantum mechanical scattering problems.
The $d$-dimensional time-independent Schr\"odinger equation for a system with unit mass is given by
\begin{equation}
  \left(-\frac{1}{2} \Delta  + V(\mathbf{x}) - E\right) \psi(\mathbf{x}) = \phi(\mathbf{x}), \quad \text{for} ~ \textbf{x} \in \mathbb{R}^d, 
\end{equation}
where $\Delta$ is the $d$-dimensional Laplacian, $V(\mathbf{x})$ is
a scalar potential, $\psi$ is the wave function and $\phi$ is the right-hand side, 
which is often related to the ground state of the system. 
Depending on the total energy $E \in \mathbb{R}$, the above system allows for scattering solutions, 
in which case the equation can be reformulated as a Helmholtz equation of the form
\begin{equation}
(-\Delta-k^2(\mathbf{x})) \, \psi(\mathbf{x}) = 2\phi(\mathbf{x}), \quad \text{for} ~ \textbf{x} \in \mathbb{R}^d, 
\end{equation}
where the spatially dependent wavenumber $k(\mathbf{x})$ is defined by
$k^2(\mathbf{x}) = 2(E-V(\mathbf{x}))$.  The experimental observations
from this type of quantum mechanical systems are typically far field
maps of the solution \cite{vanroose2006double}.  Indeed, in an
experimental setup, detectors are commonly placed at large distances
from the object compared to the size of the system.  These detectors
consequently measure the probability of particles escaping from the
system in certain directions.  In many quantum mechanical systems the
potential $V(\mathbf{x})$ is an analytical function, which suggests
analyticity of the wavenumber $k(\mathbf{x})$ in the above Helmholtz
equation.  Additionally, the potential is often a decaying function.
Hence, for these types of problems, the wavenumber naturally satisfies
all conditions for the use of the proposed complex contour method to
efficiently calculate the corresponding far field map.

In the paragraphs below, we first discuss a 2D model problem where single and
double ionization occur, corresponding to waves describing respectively a 
single particle or two particles escaping from the quantum mechanical system. 
The first leads to very localized evanescent waves that propagate along the boundaries of the 
computational domain, while the latter gives rise to waves traversing the full domain.
The corresponding 2D Schr\"odinger problem will be solved on a discretized 
numerical domain for a range of energies $E$ below and above
the double ionization energy threshold. For each energy, we extract the
single and double ionization cross sections, which correspond to probabilities 
of particles escaping from the system, with the help of an integral 
of the Green's function over the numerical box. The cross sections are calculated 
using both a traditional method, where the Helmholtz equation is solved on a standard 
ECS-bounded grid \cite{mccurdy2004solving}, and the new complex contour method, introduced in Section \ref{sec:Far}. 

The main purpose of the calculations in Sections \ref{subsec:2dsystem}
to \ref{subsec:validation2D} is to validate the complex contour method
when applied to Schr\"odinger's equation.  In Section
\ref{subsec:MG2D} the multigrid performance for solving the 2D
complex-valued scattered wave system is benchmarked.  Note, however,
that these 2D problems essentially do not require a multigrid solver,
since a direct sparse solver performs well for these relatively
small-size problems. In Section \ref{subsec:3dsystem} we illustrate the
convergence of a multigrid solver on a 3D Schr\"o\-ding\-er equation, with
energies that allow for triple ionization as well as double and
single ionization. The previous attempts to solve these problems with
the help of the complex shifted Laplacian as a preconditioner to a general 
Krylov method showed a notable deterioration in the convergence behavior in function of the total energy $E$                                                                                                
\cite{reps2010indefinite}. 
It will be shown that the new complex contour method, which allows to use 
multigrid as a solver, performs well for these problems.

Although the benchmark problems considered in this section mainly use model potentials, 
we believe that the calculations presented below are an important step towards the application 
of the method on realistic quantum mechanical systems.


\subsection{Cross sections of the 2D Schr\"odinger problem}\label{subsec:2dsystem}

Our primary aim is to validate the applicability of the new complex
contour method on the 2D Schr\"odinger equation describing a quantum
mechanical scattering problem. This problem originates from the
expansion of a 6D scattering problem in spherical harmonics, see
\cite{baertschy2001electron,vanroose2006double}, in which each
particle is expressed in terms of its spherical coordinates, 
resulting in a coupled system of 2D equations. Note that the
differential operators only appear on the diagonal blocks of the system.  
These diagonal blocks then form the two-dimensional Schr\"odinger equation
\begin{equation}\label{eq:partialwave2d}
\left(-\frac{1}{2} \Delta  + V_1(x) +  V_2(y) + V_{12}(x,y) -E\right) u(x,y) = \phi(x,y),   \quad x,y \ge 0,
\end{equation}
with boundary conditions
\begin{equation}\label{eq:partialwave2dBC}
\begin{cases} 
u(x,0)=0 \quad \text{for} \quad x\ge 0\\
u(0,y)=0 \quad \text{for} \quad y\ge 0\\
\text{outgoing} \quad \text{for} \quad x \rightarrow \infty \quad \text{or} \quad y \rightarrow \infty,
\end{cases}
\end{equation}
where $\Delta$ represents the 2D Laplacian, $V_1(x)$ and $V_2(y)$ are
the one-body potentials, $V_{12}(x,y)$ is a two-body potential and $E$
is the total energy of the system. Since the arguments $x$ and $y$ are
in fact radial coordinates in the partial wave expansion, homogeneous
Dirichlet boundary conditions are implied at the $x=0$ and $y=0$
boundaries. The potentials $V_1$, $V_2$ and $V_{12}$ are generally
analytical functions that decay as the radial coordinates $x$ and $y$
become large.

Depending on the strength of the one-body potentials $V_1$ and $V_2$,
the problem allows for so-called single ionization waves, which are localized evanescent
waves that propagate along the edges of the domain. We refer the reader to Sections \ref{subsec:spectralprops} and \ref{subsec:solutypes}
for a more detailed physical clarification. In the following, we expound on the situation 
with a strong attractive potential $V_1$ in the $x$-direction; the case with a strong $V_2$ 
potential is completely analogous. If the attraction of $V_1$ is strong enough, there exists a 
one-dimensional eigenstate $\phi_n(x)$ for every negative eigenvalue $\lambda_n < 0$, 
characterized by a one-dimensional Helmholtz equation
\begin{equation}\label{eq:boundstate}
  \left(-\frac{1}{2}\frac{d^2}{dx^2}+V_1(x)\right)\phi_n(x) = \lambda_n \, \phi_n(x),  \quad \text{for} ~ x \ge 0.
\end{equation}
Note that $\phi_n(0)=0$ and $\phi_n(x \rightarrow \infty)=0$. 

The far field maps of this system are then again Green's integrals over the solution, see
\cite{mccurdy2004solving}. Indeed, the single ionization amplitude $s_n(E)$, which represents the total number of single ionized particles, is given by
\begin{equation} \label{eq:single-ionization}
  s_n(E) = \int_\Omega  \phi_{k_n}(x) \phi_n(y) \left[\phi(x,y) -  V_{12}(x,y) u(x,y)\right]  \, dx \, dy,
\end{equation}
where $k_n = \sqrt{2(E-\lambda_n)}$, $\phi_n$ is a one-body
eigenstate, i.e.\, a solution of equation \eqref{eq:boundstate} with a corresponding eigenvalue
$\lambda_n$, and the function $\phi_{k_n}$ is a regular, normalized solution of the homogeneous Helmholtz equation
\begin{equation}\label{eq:testfunction}
\left(-\frac{1}{2}\frac{d^2}{dx^2} + V_1(x) -\frac{1}{2}k^2\right)\phi_{k} = 0,
\end{equation}
where $k=k_n$ and $\phi_{k_n}$ normalized with $1/\sqrt{k_n}$. 

Similarly, the double ionization cross section $f(k_1,k_2)$, which measures the total
number of double ionized particles, is defined by the integral
\begin{equation} \label{eq:double-ionization}
 f(k_1,k_2) =  \int_\Omega  \phi_{k_1}(x)  \phi_{k_2}(y)\left(\phi(x,y) -  V_{12}(x,y) u(x,y)\right) \, dx \, dy,  \quad x,y \ge 0,
\end{equation}
where both $\phi_{k_1}(x)$ and $\phi_{k_2}(y)$ are solutions to
\eqref{eq:testfunction}, with $k_1 = \sqrt{2E} \sin(\alpha)$ and $k_2
= \sqrt{2E}\cos(\alpha)$, respectively, and $\alpha \in [0,\pi/2]$
such that $k_1^2 + k_2^2 = 2E$.  The total double ionization cross
section is defined as the integral
\begin{equation}
  \sigma_{tot}(E) = \int_0^E \sigma(\sqrt{2\epsilon},\sqrt{2(E-\epsilon)}) \, d\epsilon,
\end{equation}
where 
\begin{equation}
  \sigma(k_1,k_2) = \frac{8\pi^2}{k_0^2} \frac{1}{k_1k_2} |f(k_1,k_2)|^2.
\end{equation}
The above integral expressions are obtained through a reorganization similar 
to the one performed in Section \ref{sec:Helmholtz}, see \eqref{eq:helmholtz_reorganized}-\eqref{eq:HHR}. 
For example, to calculate the single ionization cross section, 
equation \eqref{eq:partialwave2d} is reorganized as
\begin{equation}
  \left(-\frac{1}{2} \Delta  + V_1(x) -E\right) u(x,y) = \phi(x,y) - (V_2(y) + V_{12}(x,y)) \, u(x,y),   \quad x,y \ge 0.
\end{equation}
Since the left hand side is separable, the corresponding Green's function allows us to write 
\begin{equation}
  u(x,y)=\int_\Omega G(x,y|x',y') \left(\phi(x',y')-(V_2(y')+V_{12}(x',y')) \, u^N(x',y')\right) \, dx' \, dy'.
\end{equation}
Using the asymptotic form of the Green's function, the above ultimately results 
in integral formulation \eqref{eq:single-ionization}. The double ionization 
integral expression \eqref{eq:double-ionization} can be derived in a similar way.

\subsection{Spectral properties of the 2D Schr\"odinger
  problem}\label{subsec:spectralprops}
\begin{figure}
\begin{center}
\begin{tabular}{cc}
 \includegraphics[width=0.48\textwidth]{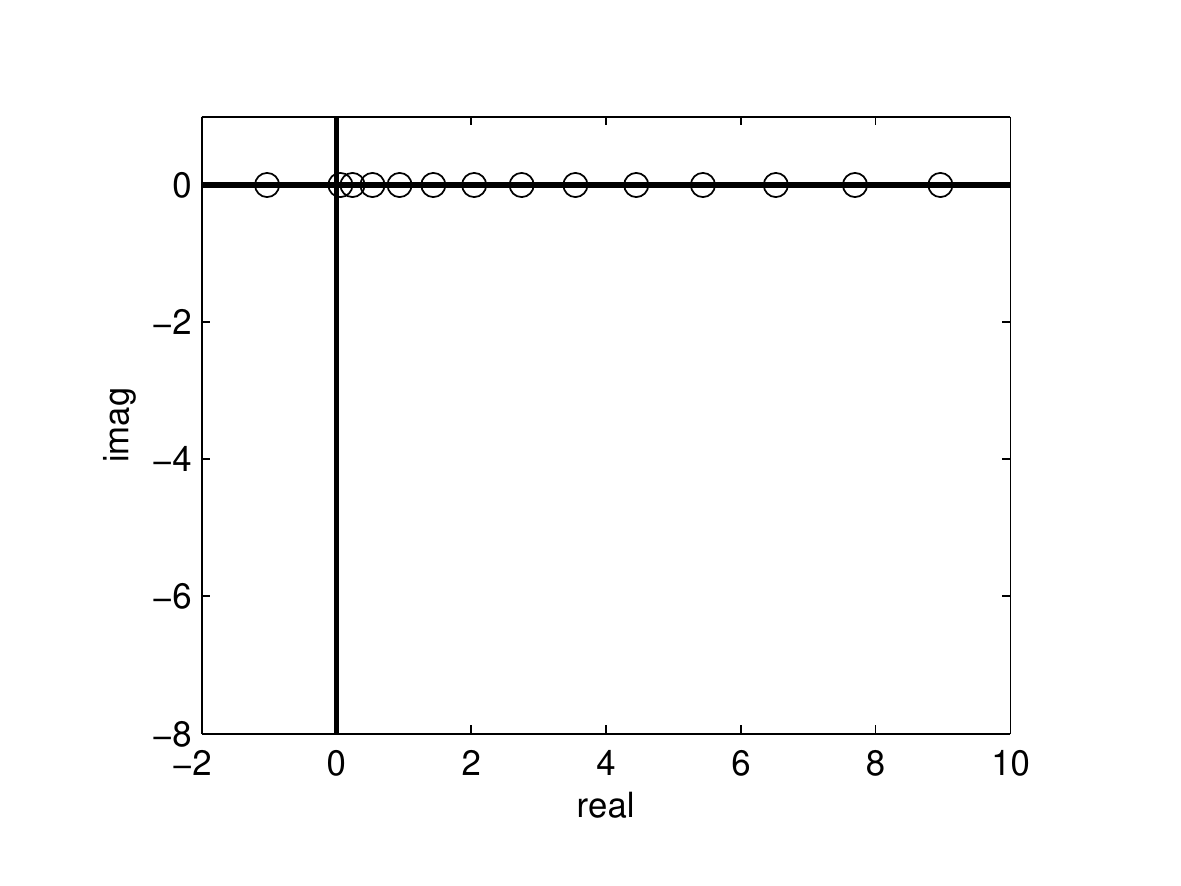} & 
 \includegraphics[width=0.48\textwidth]{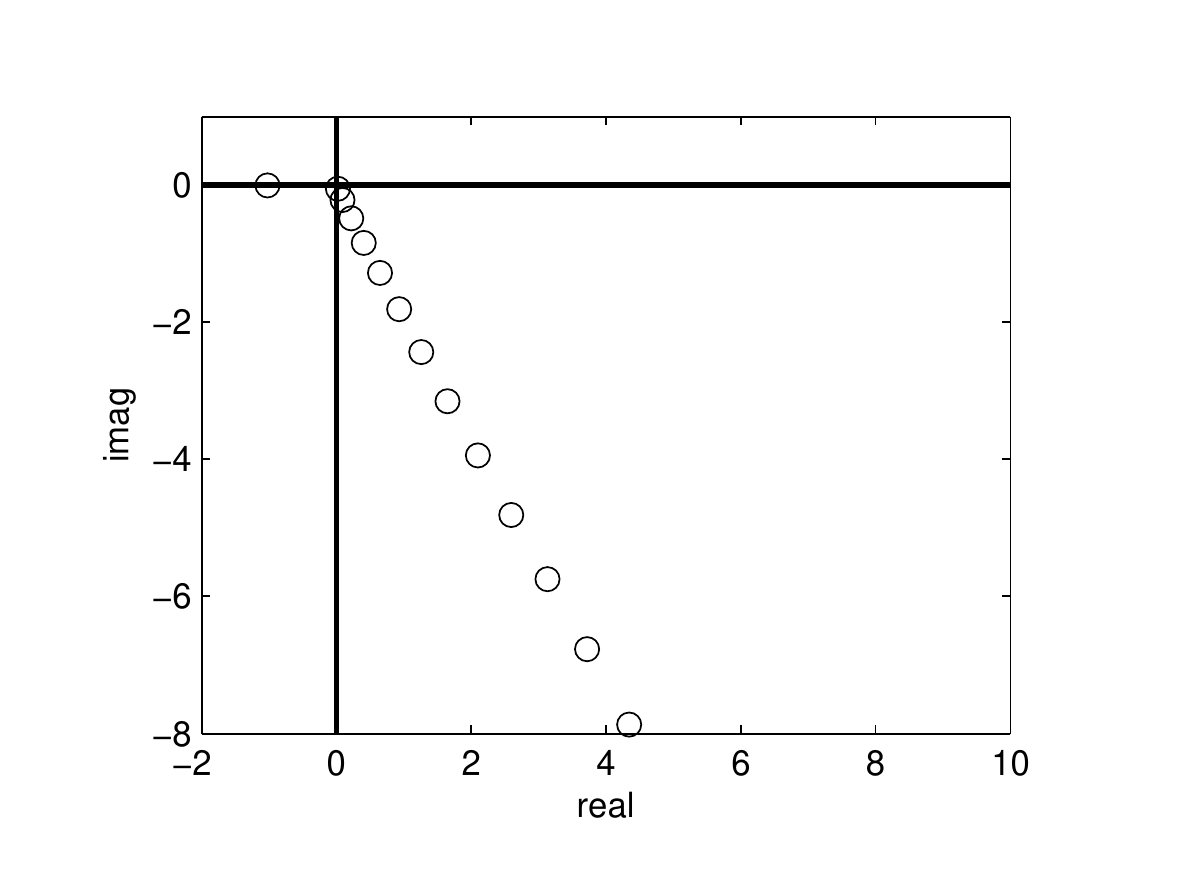}\\
 \includegraphics[width=0.48\textwidth]{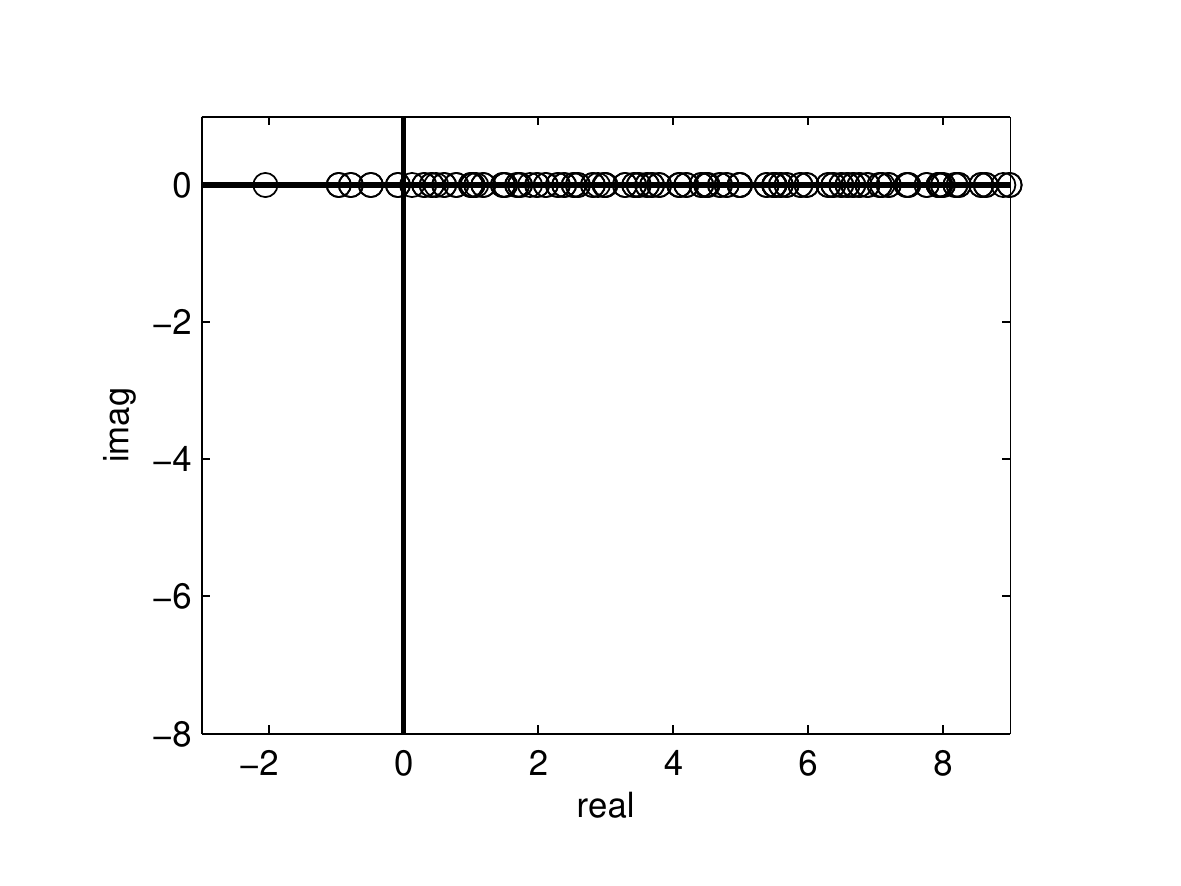} & 
 \includegraphics[width=0.48\textwidth]{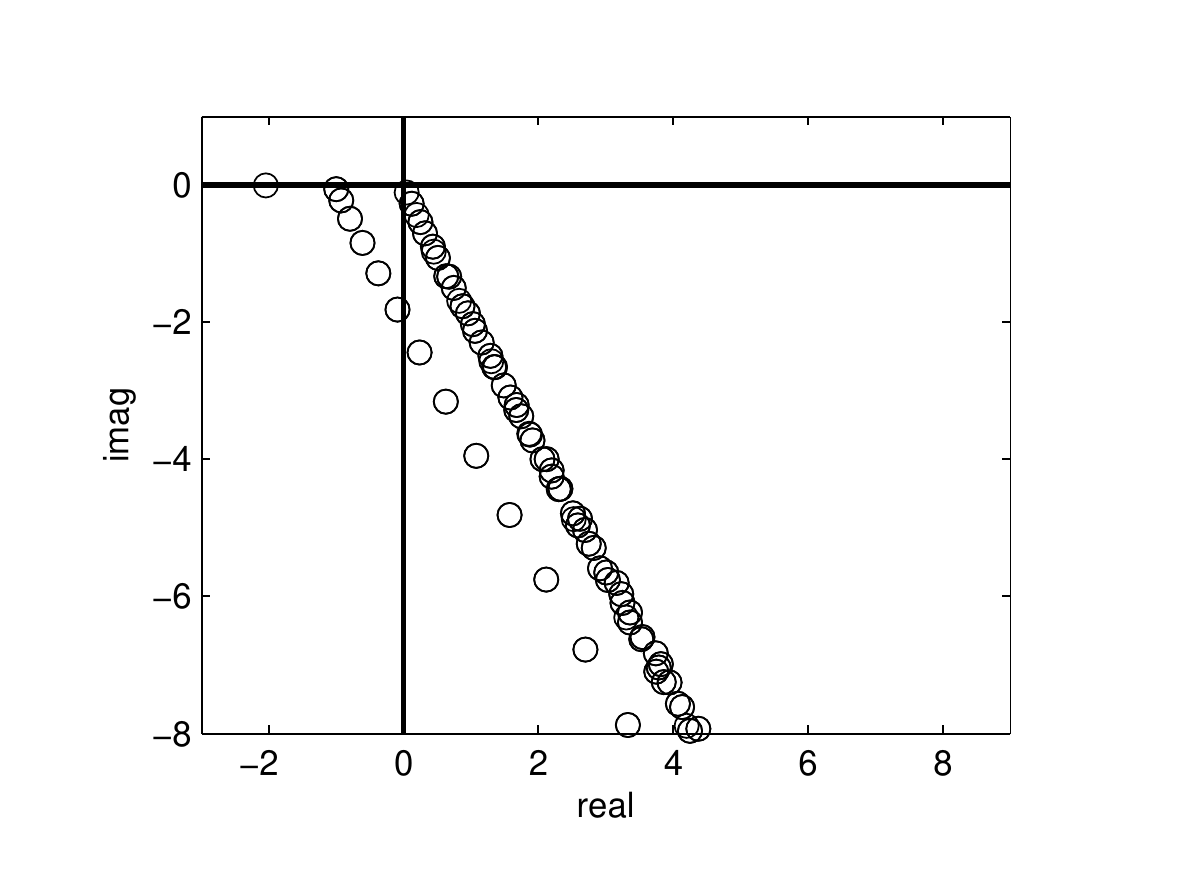}
\end{tabular}
\end{center}
\caption{Spectrum (close-up) of the discretized 1D (top) and 2D (bottom) Hamiltonian (\ref{eq:hamiltonian2d}), i.e.\ with $E = 0$. Left: spectrum corresponding to the standard discretization with real-valued grid distance $h$. Right: spectrum for a complex contour discretization with complex-valued grid distance $\tilde{h} = e^{-i\gamma} h$, where $\gamma = \pi/6$. The complex spectrum is rotated down into the complex plane over $2\gamma = \pi/3$, cf.\ \cite{reps2010indefinite}.}
\label{fig:eigenvalues}
\end{figure}

To obtain insight in the multigrid convergence for the two-dimensional Schr\"o\-ding\-er problem, 
we briefly discuss the spectral properties of the discretized
Schr\"o\-ding\-er operator. The discretized 2D Hamiltonian $H^{2d}$ corresponding to equation
\eqref{eq:partialwave2d} can be written as a sum of two Kronecker
products and a two-body potential, i.e.
\begin{equation}
  H^{\text{2d}} = H^{\text{1d}}\otimes I + I \otimes H^{\text{1d}} +  V_{12}(x,y),
\end{equation}
where $H^{\text{1d}} = -1/2 \Delta + V_i$ $(i = 1,2)$ is the one-dimensional Hamiltonian, discretized 
using finite differences. When the two-body potential $V_{12}(x,y)$ is weak relative to the one-body 
potentials, the eigenvalues of the 2D Hamiltonian can be approximated by
\begin{equation} \label{eq:hamiltonian2d}
  \lambda^{\text{2d}} \approx \lambda_i^{\text{1d}} + \lambda_j^{\text{1d}}, \qquad 1 \leq i,j, \leq n.
\end{equation}
Hence, to form a better understanding of the spectral properties of $H^{\text{2d}}$, 
let us first consider the eigenvalues of $H^{1d}$. After discretization using 
second order finite differences, the one-dimensional Hamiltonian can be written as a tridiagonal matrix,
where the stencil 
\begin{equation} 
\frac{1}{h^2}\left[
\begin{matrix} 
 -1 & 2 & -1 
\end{matrix} \right]
\end{equation}
approximates the second
derivatives, and the potential is a diagonal matrix evaluated in the
grid points. The spectrum of $H^{\text{1d}}$ closely resembles the 
spectrum of the Laplacian $(-1/2) \Delta$, however the presence of the
potential modifies the smallest eigenvalues. The resulting spectrum is 
shown on the top left panel of Figure \ref{fig:eigenvalues}, which 
presents a close-up of the eigenvalues near the origin. A single negative
eigenvalue $\lambda_0^{1d}=-1.0215$ can be observed, which is due to the attractive 
potential. The remaining spectrum consists of a series of positive 
eigenvalues located along the positive real axis.
The top right panel of Figure \ref{fig:eigenvalues} shows the eigenvalues of $H^{1d}$,
discretized along a complex-valued contour, i.e.\ the real grid rotated into the complex plane by $\gamma$.  
The grid distance used is now $\tilde{h}=h e^{i\gamma}$, which results in the following
stencil for the second derivative
\begin{equation} 
e^{-2i\gamma}\frac{1}{h^2}\left[
\begin{matrix} 
 -1 & 2 & -1 
\end{matrix} \right].
\end{equation}
This implies that the spectrum of the Laplacian is rotated down into 
the complex plane by an angle $2\gamma$. Figure \ref{fig:eigenvalues} 
shows that most of the eigenvalues of $H^{\text{1d}}$ are rotated downwards over $2\gamma$, 
with the exception of the bound state eigenvalue $\lambda_0^{1d}$, which remains on the negative 
real axis. The changes of  spectrum as a result of a rotation are well known in the physics literature, 
see for example \cite{moiseyev1998quantum}. 

We consequently turn to the two-dimensional problem setting, where the
eigenvalues of the Hamiltonian $H^{\text{2d}}$ are approximately
sums of the one-dimensional operator eigenvalues $H^{\text{1d}}
\otimes I + I \otimes H^{\text{1d}}$.  The resulting eigenvalues are
shown on the bottom two panels of Figure~\ref{fig:eigenvalues}. Again,
the eigenvalues of the 2D Hamiltonian are rotated down in the complex
plane when the system is discretized along a complex-valued
contour. In 2D, an isolated eigenvalue appears around $\lambda_0^{1d}
+ \lambda_0^{1d} = -2.043$, and two series of eigenvalues emerge from
the real axis: a first branch of eigenvalues starting at $-1.012$,
which originates from the sum of the negative eigenvalue
$\lambda_0^{1d}$ of the first 1D Hamiltonian combined with all the
positive eigenvalues of the second 1D Hamiltonian; and a second series
of eigenvalues starting at the origin, originating from the sums of
the positive eigenvalues of both one-dimensional Hamiltonians.

\subsection{Solution types: single and double
  ionization}\label{subsec:solutypes}
The Schr\"odinger equation \eqref{eq:partialwave2d} can be written
shortly as $(H-E\,I)u=\phi$, where $H = (-1/2)\Delta + V_1 + V_2 +
V_{12}$ is the Hamiltonian and the scalar $E$ is the total energy of
the system.  Depending on this energy $E$, the Schr\"odinger system
has different types of solutions.  In this section we briefly expound
on the physical interpretation of these solution types, using a model
problem example.  This section may prove less interesting to readers
who are primarily interested in the computational aspects of the
solution, and as such can be skipped at will.

When the energy $E$ is larger than the smallest (negative) eigenvalue of 
the Hamiltonian, i.e.\ $\lambda_0 < E$, one-body 
eigenstate solutions to \eqref{eq:boundstate} arise. 
These eigenstates can be combined into separable waves of the form
\begin{equation}\label{eq:single_ionization}
  u_n(x,y) =  \phi_n(x) \exp(ik_n y),
\end{equation}
with $k_n = \sqrt{2(E-\lambda_n)}$, that form a solution to the
Schr\"odinger problem \eqref{eq:partialwave2d} in the region where $x$
is small and $y \rightarrow \infty$.  Indeed, when $y$ is large, the
potentials $V_2$ and $V_{12}$ are negligibly small and the resulting
Schr\"odinger equation becomes separable in variables, where
$\phi_n(x)$ and $\exp(ik_n y)$ are the solutions of the separated
operators respectively.  An analogous argument holds for the case when
$V_2$ dominates and $V_1$ and $V_{12}$ are negligibly small, 
yielding evanescent waves of the form
\begin{equation} \label{eq:single_ionization2}
  u_n(x,y) =  \phi_n(y)\exp(ik_n x).
\end{equation}
These separable waves are solutions of the Schr\"odinger system for 
$x\rightarrow \infty$ and $y$ small, and can be derived similarly to 
\eqref{eq:single_ionization}. 

\begin{figure}
\begin{center}
\begin{tabular}{cc}
 {\scriptsize Single ionization ~~} & {\scriptsize Double ionization ~~} \vspace{-0.1cm}\\
 \includegraphics[width=0.48\textwidth]{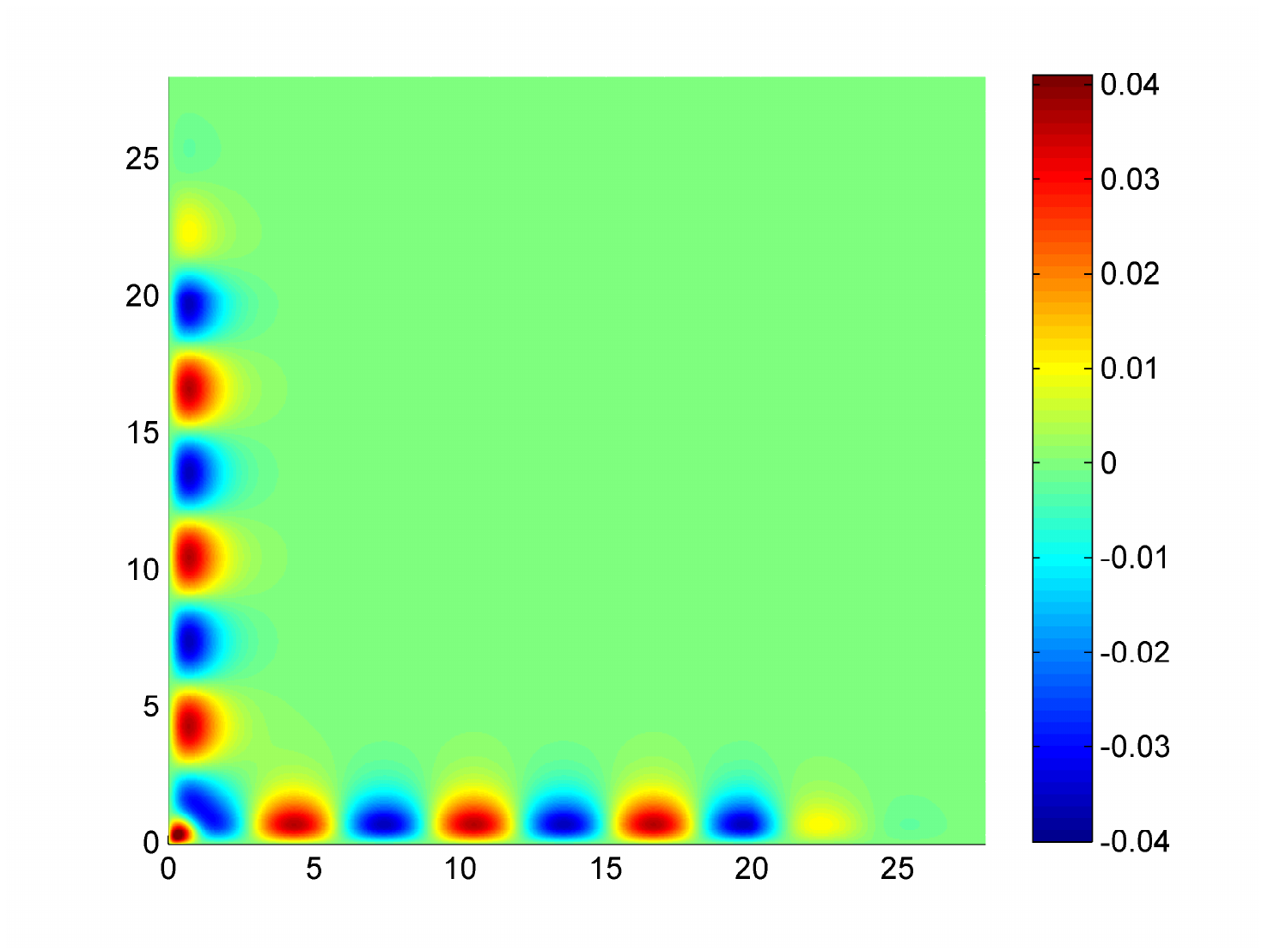} & 
 \includegraphics[width=0.48\textwidth]{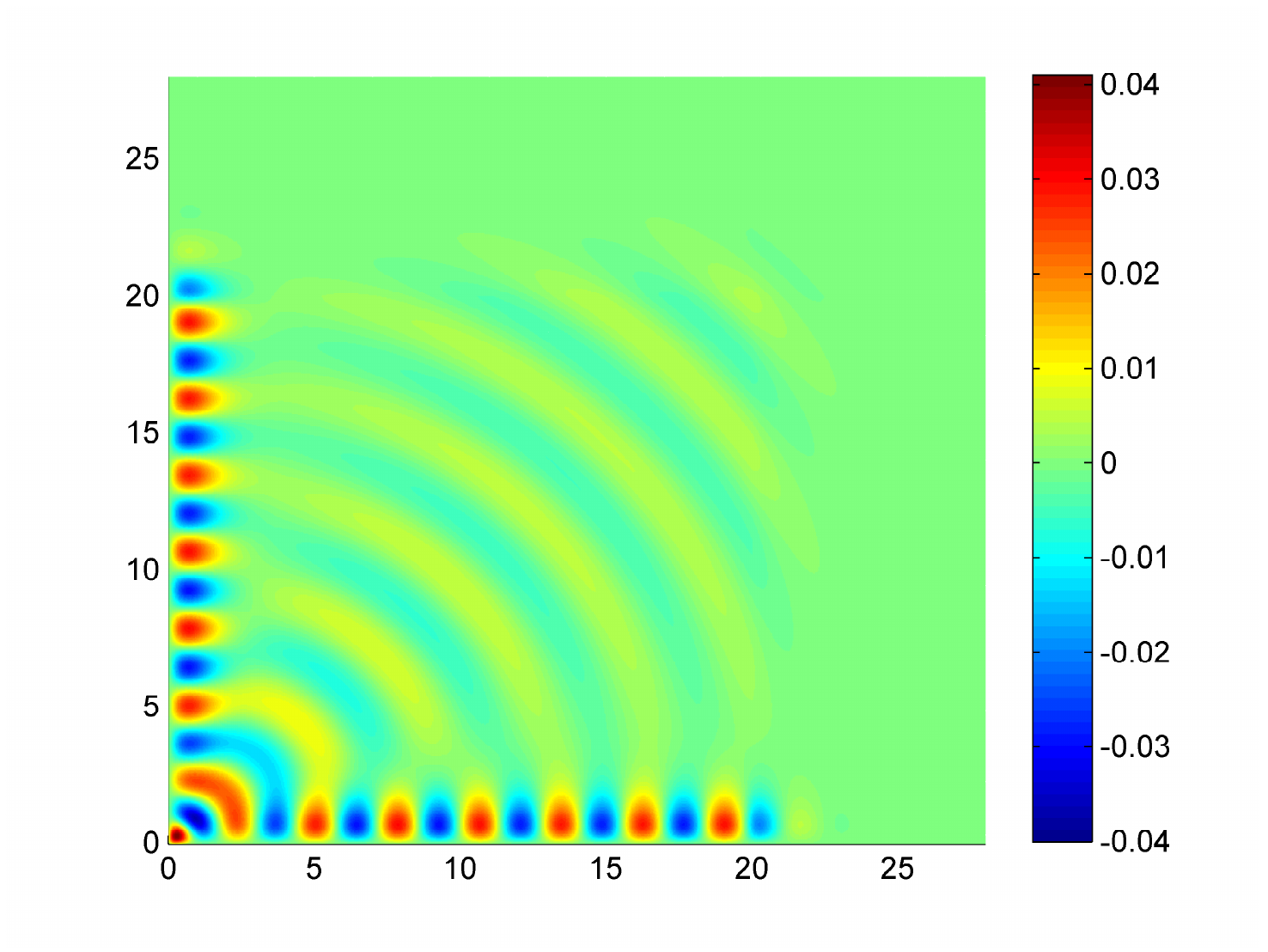}
\end{tabular}
\end{center}
\caption{Scattered wave solutions $u(x,y)$ to model problem (\ref{eq:partialwave2d}). 
	Model specifications: see accompanying text. Left: solution for energy $E=-0.5$ where only single
  ionization occurs. Right: solution for energy $E=1.5$ where both double and single
  ionization occur. Single ionization waves are localized solely along the
  edges of the domain (left), while double ionization waves appear both
  along the edges and in the middle of the domain (right).}
\label{fig:example1_solutions}
\end{figure}

Note that we can associate such a separable wave with each eigenstate 
$\phi_n$ of equation \eqref{eq:boundstate} that corresponds to a negative
eigenvalue of the Hamiltonian. These localized waves are called 
single ionization waves in the physics literature, since they
correspond to a quantum mechanical system in which a single particle is
ionized. Single ionization waves are present in the solution as soon as the
energy $E$ is above the $\lambda_0$ threshold. When there is a second
eigenstate with negative energy, say $\phi_1$ with $\lambda_1$, an
additional single ionization wave appears in the problem as soon as
$E> \lambda_1$. We refer to the specialized literature for a detailed
discussion of the ionization process, see
\cite{friedrich1991theoretical,newton2002scattering}.

The left panel of Figure~\ref{fig:example1_solutions} shows the solution 
to \eqref{eq:partialwave2d} for a total energy $E=-0.5$. 
The model problem under consideration fits
equation \eqref{eq:partialwave2d}, with a right-hand side given by 
$\phi(x,y) = \exp(-3(x+y)^2)$. The one-body potentials are defined by 
$V_1(x) = -4.5\exp(-x^2)$ and $V_2(y) = -4.5\exp(-y^2)$, yielding an eigenstate of
equation \eqref{eq:boundstate} with energy $\lambda_0= -1.0215$.
The two-body potential equals $V_{12}(x,y)= {2\exp(-(x+y)^2)}$.
The equation is discretized on a $[0,20]^2$ domain using 500 grid points 
in every spatial dimension.
An ECS absorbing layer consisting of an additional 250 grid points 
damps the outgoing waves along the right and top edges of the domain.
Note from Figure \ref{fig:example1_solutions} how the single ionization 
eigenstate solutions given by \eqref{eq:single_ionization}-\eqref{eq:single_ionization2} 
appear along the edges of the domain. 

When the total energy $E>0$, additional asymptotic scattering
solutions to the Schr\"odinger equation appear.  If all potentials are
asymptotically zero, equation \eqref{eq:partialwave2d} boils down to a
Helmholtz equation with wavenumber $k=\sqrt{2E}$ for $x\rightarrow
\infty$ and $y \rightarrow \infty$. The resulting waves are known as
double ionization waves and physically correspond to the
simultaneous ejection of two particles from the quantum mechanical
system. In the far field, double ionization waves behave as 
$e^{i \sqrt{2E} \sqrt{x^2+y^2}}$, with an angle-dependent prefactor.  At
the same time, it is still possible to have single
ionization for positive energies $E$, since $E > \lambda_0$.  The right panel of Figure
\ref{fig:example1_solutions} shows the solution for $E = 1.5$, clearly displaying
the double ionization waves in the middle of the domain and single
ionization along the edges, where the two solutions coexist.
Additionally, one observes that the single ionization waves oscillate
faster in the $x$- or $y$-direction than a free wave with wavenumber
$k=\sqrt{2E}$, since $k_n \geq k$.

\begin{figure}
\begin{center}
\begin{tabular}{cc}
\includegraphics[width=0.48\textwidth]{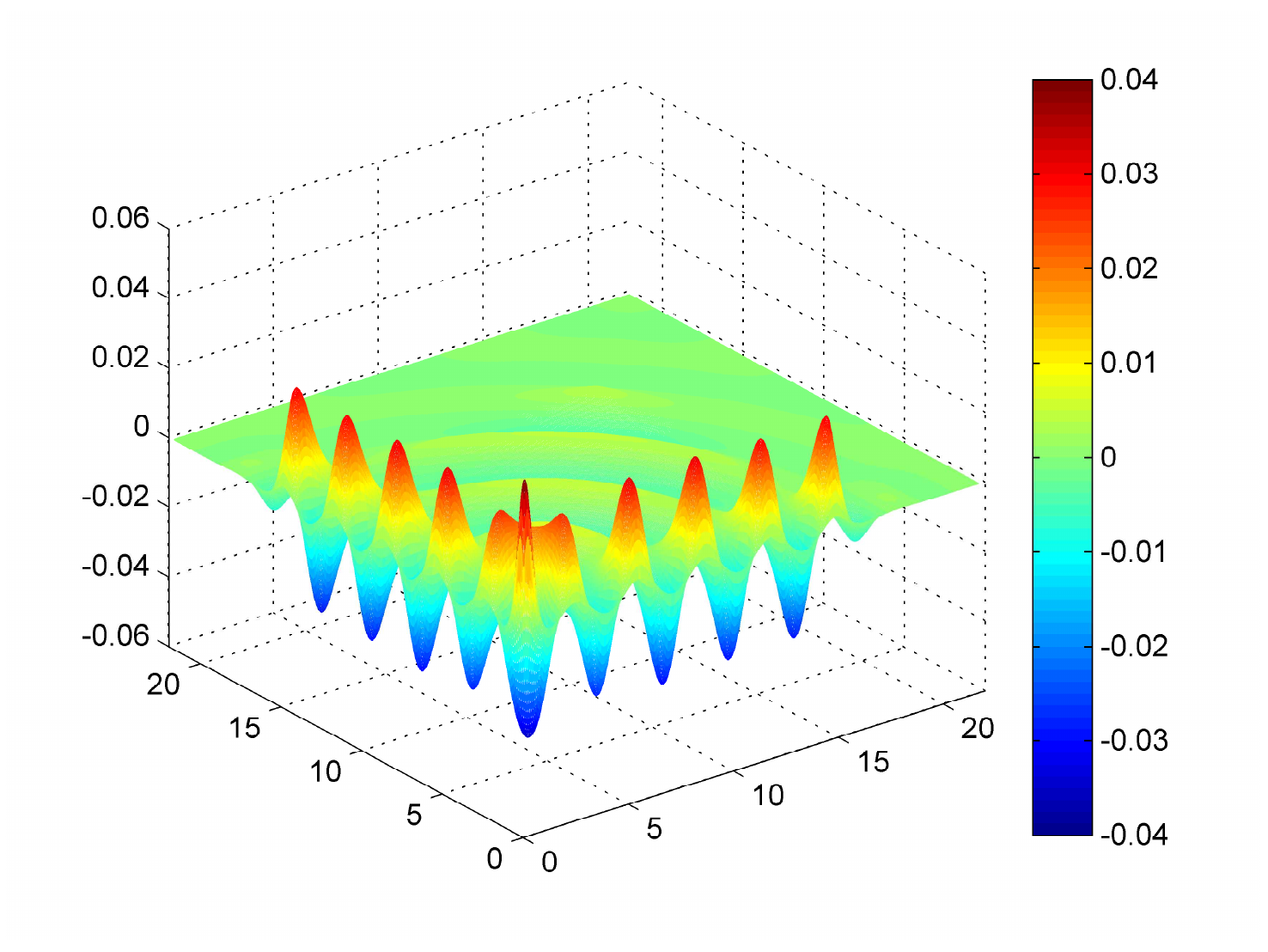} & 
\includegraphics[width=0.48\textwidth]{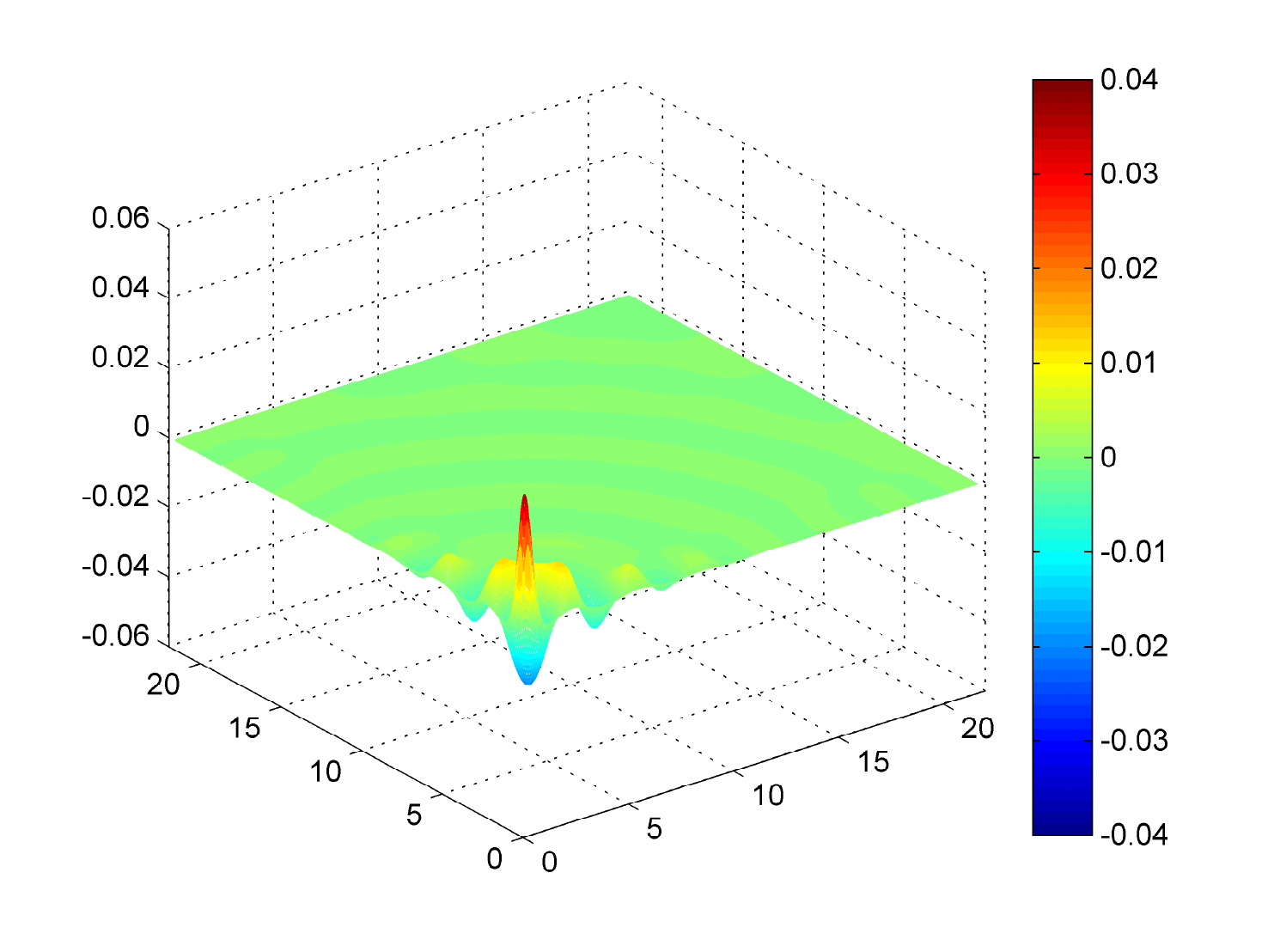}
\end{tabular}
\end{center}
\caption{Scattered wave solutions $u(x,y)$ to model problem (\ref{eq:partialwave2d}) for $E = 1$. 
  Left: solution on a real grid with ECS absorbing 
  boundary layer $(\theta = \pi/7 \approx 25^\circ)$. Right: solution 
  on a straight complex scaled contour $(\gamma \approx 8.5^\circ)$, resulting in a damped scattered wave solution.}
\label{fig:schroedinger_solution}
\end{figure}

\subsection{Complex contour method for the 2D Schr\"odinger problem} \label{subsec:validation2D}

In physical experiments, the total number of single ionized or double
ionized particles is typically observed for a range of energy levels
$E$ using advanced detectors.  These observations are made far away
from the object and effectively measure the far field amplitudes of
the solutions. The outcome of this type of experiments can be
predicted by calculating the single \eqref{eq:single-ionization} and
double ionization \eqref{eq:double-ionization} cross sections, using
the numerical solution of equation (\ref{eq:partialwave2d}), see
\cite{mccurdy2004solving}.  Indeed, in order to calculate these cross
sections, the numerical solution $u^N(x,y)$ to
\eqref{eq:partialwave2d} is required, which is generally hard to
obtain, especially in higher spatial dimensions.  However, since the
potentials $V_1$, $V_2$ and $V_{12}$ are analytical functions, the
integrals for single and double ionization can be calculated along a
complex contour rather than the classical real domain, in analogy to
the discussion in Section \ref{sec:Far}.  Consequently, one requires
the scattering solution of equation \eqref{eq:partialwave2d} on a
complex contour, which is much easier to compute. 

In the following, we calculate the single and double ionization cross
sections for a number of energies $E$ between $-1$ and $3$ using both
the classical real-valued discretization and the complex contour
approach proposed in this paper. The corresponding 2D scattering
problems (\ref{eq:partialwave2d}) are solved on a numerical domain
$\Omega = [0,15]^2$ covered by a finite difference grid consisting of
300 grid points in every spatial dimension.  Additionally, an ECS
absorbing boundary layer starts at $x = 15$ and $y = 15$,
respectively, and implements the outgoing boundary conditions,
adding an additional 150 grid points in every spatial dimension. The
ECS angle is $\theta = \pi/7 \approx 25.7^\circ$.  For the complex contour method, a
complex scaled grid with an overall complex rotation angle $\gamma
\approx 8.5^\circ$ is used.  Solutions $u^N(x,y)$ to
\eqref{eq:partialwave2d} for a total energy $E = 1$ on both the
classical real-valued grid and the complex contour are presented on
Figure \ref{fig:schroedinger_solution}.  Note how the solution is
damped when evaluated along the complex contour.

\begin{figure}
\begin{center}
 \begin{tikzpicture}
\begin{axis}[
        xlabel=$E$,
        ylabel=Cross section,
        ymax=0.0045]

                \addplot[smooth,black,dashed] plot coordinates{
(-1.0215007e+00   ,  0 		   )
(-9.2150069e-01   ,  6.3858252e-04)
(-8.2150069e-01   ,  9.2467986e-04)
(-7.2150069e-01   ,  1.0749516e-03)
(-6.2150069e-01   ,  1.1890037e-03)
(-5.2150069e-01   ,  1.2934345e-03)
(-4.2150069e-01   ,  1.3972521e-03)
(-3.2150069e-01   ,  1.4765268e-03)
(-2.2150069e-01   ,  1.5438088e-03)
(-1.2150069e-01   ,  1.6045563e-03)
(-2.1500689e-02   ,  1.6558229e-03)
( 7.8499311e-02   ,  1.6941891e-03)
( 1.7849931e-01   ,  1.7290531e-03)
( 2.7849931e-01   ,  1.7594438e-03)
( 3.7849931e-01   ,  1.7788039e-03)
( 4.7849931e-01   ,  1.7946639e-03)
( 5.7849931e-01   ,  1.8132298e-03)
( 6.7849931e-01   ,  1.8252779e-03)
( 7.7849931e-01   ,  1.8288557e-03)
( 8.7849931e-01   ,  1.8339955e-03)
( 9.7849931e-01   ,  1.8431227e-03)
( 1.0784993e+00   ,  1.8474371e-03)
( 1.1784993e+00   ,  1.8438096e-03)
( 1.2784993e+00   ,  1.8406952e-03)
( 1.3784993e+00   ,  1.8437351e-03)
( 1.4784993e+00   ,  1.8469109e-03)
( 1.5784993e+00   ,  1.8426860e-03)
( 1.6784993e+00   ,  1.8337866e-03)
( 1.7784993e+00   ,  1.8289090e-03)
( 1.8784993e+00   ,  1.8302016e-03)
( 1.9784993e+00   ,  1.8307270e-03)
( 2.0784993e+00   ,  1.8241433e-03)
( 2.1784993e+00   ,  1.8130237e-03)
( 2.2784993e+00   ,  1.8052771e-03)
( 2.3784993e+00   ,  1.8042338e-03)
( 2.4784993e+00   ,  1.8049649e-03)
( 2.5784993e+00   ,  1.8003346e-03)
( 2.6784993e+00   ,  1.7892010e-03)
( 2.7784993e+00   ,  1.7774952e-03)
( 2.8784993e+00   ,  1.7714823e-03)
( 2.9784993e+00   ,  1.7711985e-03)};
\addlegendentry{Single ionization (real)}  

\addplot[mark=triangle*,black,only marks] plot coordinates{
(-1.0215007e+00  ,            0 )
(-8.2150069e-01  ,  9.0468611e-04 )
(-6.2150069e-01  ,  1.2154105e-03 )
(-4.2150069e-01  ,  1.4230472e-03 )
(-2.2150069e-01  ,  1.5720915e-03 )
(-2.1500689e-02  ,  1.6830472e-03 )
( 1.7849931e-01  ,  1.7563615e-03 )
( 3.7849931e-01  ,  1.8025587e-03 )
( 5.7849931e-01  ,  1.8385737e-03 )
( 7.7849931e-01  ,  1.8504130e-03 )
( 9.7849931e-01  ,  1.8668001e-03 )
( 1.1784993e+00  ,  1.8628721e-03 )
( 1.3784993e+00  ,  1.8659478e-03 )
( 1.5784993e+00  ,  1.8603197e-03 )
( 1.7784993e+00  ,  1.8479752e-03 )
( 1.9784993e+00  ,  1.8494916e-03 )
( 2.1784993e+00  ,  1.8276939e-03 )
( 2.3784993e+00  ,  1.8228369e-03 )
( 2.5784993e+00  ,  1.8154194e-03 )
( 2.7784993e+00  ,  1.7906321e-03 )
( 2.9784993e+00  ,  1.7882720e-03 )
};
\addlegendentry{Single ionization (complex)}

\addplot[smooth,black,dotted] plot coordinates{
(-1.0215007e+00 ,  0 	          ) 
(-9.2150069e-01 ,  0            )
(-8.2150069e-01 ,  0		  )
(-7.2150069e-01 ,  0		  )
(-6.2150069e-01 ,  0 	          )
(-5.2150069e-01 ,  0 	          )
(-4.2150069e-01 ,  0 	          )
(-3.2150069e-01 ,  0		  )
(-2.2150069e-01 ,  0		  )
(-1.2150069e-01 ,  0    	  )
(-2.1500689e-02 ,  0		  )
( 7.8499311e-02 ,  1.1050837e-06)
( 1.7849931e-01 ,  5.3082119e-06)
( 2.7849931e-01 ,  1.2525175e-05)
( 3.7849931e-01 ,  2.2189419e-05)
( 4.7849931e-01 ,  3.4830544e-05)
( 5.7849931e-01 ,  4.9538427e-05)
( 6.7849931e-01 ,  6.5792520e-05)
( 7.7849931e-01 ,  8.3338714e-05)
( 8.7849931e-01 ,  1.0254151e-04)
( 9.7849931e-01 ,  1.2329677e-04)
( 1.0784993e+00 ,  1.4476612e-04)
( 1.1784993e+00 ,  1.6650046e-04)
( 1.2784993e+00 ,  1.8876351e-04)
( 1.3784993e+00 ,  2.1186457e-04)
( 1.4784993e+00 ,  2.3579776e-04)
( 1.5784993e+00 ,  2.6056731e-04)
( 1.6784993e+00 ,  2.8630097e-04)
( 1.7784993e+00 ,  3.1297130e-04)
( 1.8784993e+00 ,  3.4040841e-04)
( 1.9784993e+00 ,  3.6840341e-04)
( 2.0784993e+00 ,  3.9650813e-04)
( 2.1784993e+00 ,  4.2401240e-04)
( 2.2784993e+00 ,  4.5048901e-04)
( 2.3784993e+00 ,  4.7617820e-04)
( 2.4784993e+00 ,  5.0158394e-04)
( 2.5784993e+00 ,  5.2680940e-04)
( 2.6784993e+00 ,  5.5148539e-04)
( 2.7784993e+00 ,  5.7534142e-04)
( 2.8784993e+00 ,  5.9867209e-04)
( 2.9784993e+00 ,  6.2216583e-04)
};
\addlegendentry{Double ionization (real)}

\addplot[mark=*,color=black,only marks] plot coordinates{
( -1.0215007e+00  ,  0  )
( -8.2150069e-01  ,  0  )
( -6.2150069e-01  ,  0  )
( -4.2150069e-01  ,  0  )
( -2.2150069e-01  ,  0  )
( -2.1500689e-02  ,  0  )
(  1.7849931e-01  ,   5.2937397e-06  )
(  3.7849931e-01  ,   2.2421628e-05  )
(  5.7849931e-01  ,   4.8932351e-05  )
(  7.7849931e-01  ,   8.2720542e-05  )
(  9.7849931e-01  ,   1.2196342e-04  )
(  1.1784993e+00  ,   1.6624666e-04  )
(  1.3784993e+00  ,   2.1447172e-04  )
(  1.5784993e+00  ,   2.6423830e-04  )
(  1.7784993e+00  ,   3.1460554e-04  )
(  1.9784993e+00  ,   3.6522422e-04  )
(  2.1784993e+00  ,   4.1697918e-04  )
(  2.3784993e+00  ,   4.6947186e-04  )
(  2.5784993e+00  ,   5.2275418e-04  )
(  2.7784993e+00  ,   5.7544498e-04  )
(  2.9784993e+00  ,   6.2731293e-04  )
}
;
\addlegendentry{Double ionization (complex)}

        \addplot[smooth,black] plot coordinates{
(  -1.0215007e+00,   1.1366667e-05) 
(  -9.2150069e-01,   6.9686220e-04) 
(  -8.2150069e-01,   9.0293558e-04) 
(  -7.2150069e-01,   1.0658949e-03) 
(  -6.2150069e-01,   1.1920408e-03) 
(  -5.2150069e-01,   1.2972055e-03) 
(  -4.2150069e-01,   1.3965877e-03) 
(  -3.2150069e-01,   1.4774532e-03) 
(  -2.2150069e-01,   1.5455682e-03) 
(  -1.2150069e-01,   1.6060515e-03) 
(  -2.1500689e-02,   1.6569185e-03) 
(   7.8499311e-02,   1.6981357e-03) 
(   1.7849931e-01,   1.7370135e-03) 
(   2.7849931e-01,   1.7737332e-03) 
(   3.7849931e-01,   1.8051452e-03) 
(   4.7849931e-01,   1.8346068e-03) 
(   5.7849931e-01,   1.8658103e-03) 
(   6.7849931e-01,   1.8946983e-03) 
(   7.7849931e-01,   1.9193756e-03) 
(   8.7849931e-01,   1.9444192e-03) 
(   9.7849931e-01,   1.9721788e-03) 
(   1.0784993e+00,   1.9987841e-03) 
(   1.1784993e+00,   2.0213599e-03) 
(   1.2784993e+00,   2.0432017e-03) 
(   1.3784993e+00,   2.0682682e-03) 
(   1.4784993e+00,   2.0948876e-03) 
(   1.5784993e+00,   2.1185445e-03) 
(   1.6784993e+00,   2.1387604e-03) 
(   1.7784993e+00,   2.1596717e-03) 
(   1.8784993e+00,   2.1841138e-03) 
(   1.9784993e+00,   2.2098808e-03) 
(   2.0784993e+00,   2.2327429e-03) 
(   2.1784993e+00,   2.2518907e-03) 
(   2.2784993e+00,   2.2708371e-03) 
(   2.3784993e+00,   2.2930331e-03) 
(   2.4784993e+00,   2.3178593e-03) 
(   2.5784993e+00,   2.3414250e-03) 
(   2.6784993e+00,   2.3609632e-03) 
(   2.7784993e+00,   2.3778385e-03) 
(   2.8784993e+00,   2.3959541e-03) 
(   2.9784993e+00,   2.4176249e-03) 
};
\addlegendentry{Total cross section}

\end{axis}
\end{tikzpicture}
\end{center}
\caption{Comparison of the single and double ionization total cross
  sections calculated using the scattered wave solution $u^N$ of \eqref{eq:partialwave2d} calculated on (a) a traditional 
  real-valued ECS grid with $\theta = \pi/7$ and (b) a full complex contour with $\gamma = 8.5^\circ$. 
  The energy range starts at the single ionization threshold $E=-1$, corresponding to a strictly positive cross section.
  Double ionization occurs for energy levels $E>0$. \label{fig:single_and_double} }
\end{figure}
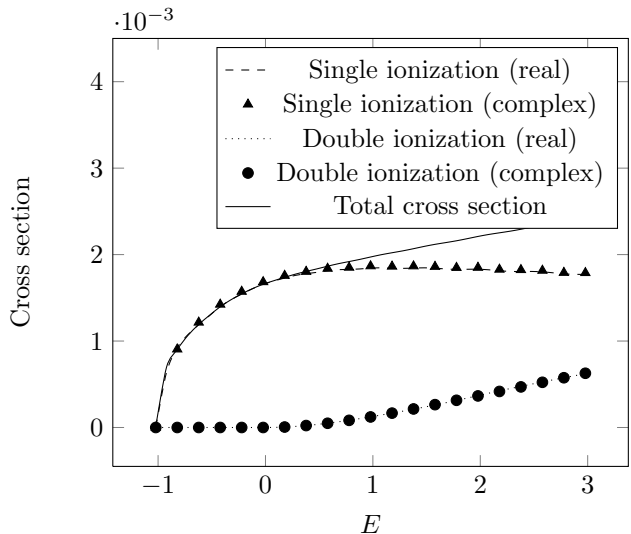

Figure \ref{fig:single_and_double} shows the rate of single and double
ionization as a function of the total energy $E$. The dashed and
dotted lines represent the single and double ionization amplitudes
calculated using the traditional real-valued method with ECS absorbing
boundary conditions \cite{mccurdy2004solving}. The solid line is the
total cross section, i.e.\, the sum of single and double ionization, and is
calculated using the optical theorem, see
\cite{newton2002scattering}. One observes that single ionization
occurs starting from $E>-1.0215$. Double ionization only occurs when
$E>0$, and comprises only a fraction of the single ionization cross
section (cf.\ Figure \ref{fig:example1_solutions}).  Note how the
energy of the single ionized bound states rises as the total energy
grows, and remains present even when $E>0$.  Results obtained using
the complex contour approach are indicated by the $\blacktriangle$ and
$\bullet$ symbols on Figure \ref{fig:single_and_double}.  In this
case, the Schr\"odinger equation \eqref{eq:partialwave2d} is first
solved on a complex contour, yielding a damped solution as shown by
Figure \ref{fig:schroedinger_solution} (right panel), followed by the
calculation of the integrals \eqref{eq:single-ionization} and
\eqref{eq:double-ionization} along this complex contour. Identical
results are obtained by both calculation methods, thus validating the
applicability of the complex contour approach on Schr\"odinger-type
problems.

\subsection{Multigrid performance on the 2D Schr\"odinger problem} \label{subsec:MG2D}

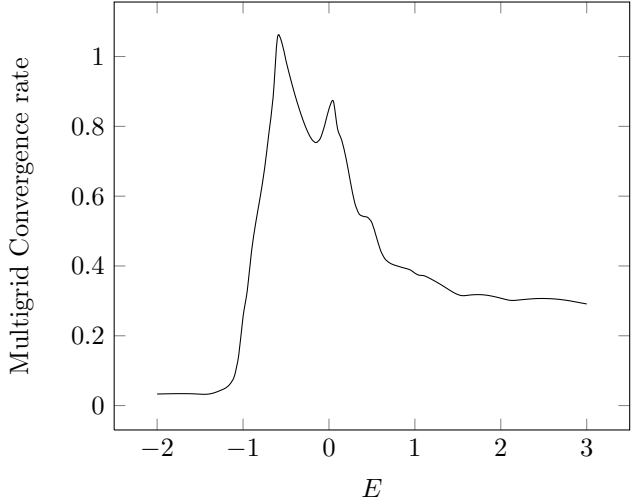
\begin{figure}
\begin{center}
 \begin{tikzpicture}
\begin{axis}[
        xlabel=$E$,
        ylabel=Multigrid Convergence rate]
        \addplot[smooth,black] plot coordinates{
(   -2.0000  ,     0.0332 )
(   -1.9500  ,     0.0334 )
(   -1.9000  ,     0.0337 )
(   -1.8500  ,     0.0339 )
(   -1.8000  ,     0.0341 )
(   -1.7500  ,     0.0342 )
(   -1.7000  ,     0.0342 )
(   -1.6500  ,     0.0341 )
(   -1.6000  ,     0.0339 )
(   -1.5500  ,     0.0334 )
(   -1.5000  ,     0.0327 )
(   -1.4500  ,     0.0323 )
(   -1.4000  ,     0.0330 )
(   -1.3500  ,     0.0356 )
(   -1.3000  ,     0.0397 )
(   -1.2500  ,     0.0449 )
(   -1.2000  ,     0.0508 )
(   -1.1500  ,     0.0622 )
(   -1.1000  ,     0.0849 )
(   -1.0500  ,     0.1448 )
(   -1.0000  ,     0.2546 )
(   -0.9500  ,     0.3308 )
(   -0.9000  ,     0.4448 )
(   -0.8500  ,     0.5273 )
(   -0.8000  ,     0.5982 )
(   -0.7500  ,     0.6778 )
(   -0.7000  ,     0.7789 )
(   -0.6500  ,     0.8836 )
(   -0.6000  ,     1.0540 )
(   -0.5500  ,     1.0413 )
(   -0.5000  ,     0.9877 )
(   -0.4500  ,     0.9389 )
(   -0.4000  ,     0.8937 )
(   -0.3500  ,     0.8529 )
(   -0.3000  ,     0.8165 )
(   -0.2500  ,     0.7856 )
(   -0.2000  ,     0.7627 )
(   -0.1500  ,     0.7535 )
(   -0.1000  ,     0.7658 )
(   -0.0500  ,     0.8030 )
(         0  ,     0.8500 )
(    0.0500  ,     0.8726 )
(    0.1000  ,     0.7931 )
(    0.1500  ,     0.7589 )
(    0.2000  ,     0.7070 )
(    0.2500  ,     0.6406 )
(    0.3000  ,     0.5808 )
(    0.3500  ,     0.5494 )
(    0.4000  ,     0.5421 )
(    0.4500  ,     0.5391 )
(    0.5000  ,     0.5238 )
(    0.5500  ,     0.4847 )
(    0.6000  ,     0.4449 )
(    0.6500  ,     0.4210 )
(    0.7000  ,     0.4097 )
(    0.7500  ,     0.4036 )
(    0.8000  ,     0.3997 )
(    0.8500  ,     0.3958 )
(    0.9000  ,     0.3925 )
(    0.9500  ,     0.3883 )
(    1.0000  ,     0.3794 )
(    1.0500  ,     0.3733   )
(    1.1000  ,     0.3726   )
(    1.1500  ,     0.3674   )
(    1.2000  ,     0.3610   )
(    1.2500  ,     0.3542   )
(    1.3000  ,     0.3469   )
(    1.3500  ,     0.3392   )
(    1.4000  ,     0.3314   )
(    1.4500  ,     0.3240   )
(    1.5000  ,     0.3181   )
(    1.5500  ,     0.3149   )
(    1.6000  ,     0.3155   )
(    1.6500  ,     0.3170   )
(    1.7000  ,     0.3177   )
(    1.7500  ,     0.3179   )
(    1.8000  ,     0.3172   )
(    1.8500  ,     0.3157   )
(    1.9000  ,     0.3134   )
(    1.9500  ,     0.3107   )
(    2.0000  ,     0.3077   )
(    2.0500  ,     0.3045   )
(    2.1000  ,     0.3018   )
(    2.1500  ,     0.3012   )
(    2.2000  ,     0.3023   )
(    2.2500  ,     0.3035   )
(    2.3000  ,     0.3046   )
(    2.3500  ,     0.3056   )
(    2.4000  ,     0.3064   )
(    2.4500  ,     0.3068   )
(    2.5000  ,     0.3068   )
(    2.5500  ,     0.3065   )
(    2.6000  ,     0.3059   )
(    2.6500  ,     0.3050   )
(    2.7000  ,     0.3037   )
(    2.7500  ,     0.3021   )
(    2.8000  ,     0.3002   )
(    2.8500  ,     0.2980   )
(    2.9000  ,     0.2957   )
(    2.9500  ,     0.2934   )
(    3.0000  ,     0.2912   )
};                         
\end{axis}
\end{tikzpicture}
\end{center}
\caption{2D Schr\"odinger problem \eqref{eq:partialwave2d} for a total energy range $E \in [-2,3]$ solved on a full complex grid with $\gamma \approx 8.5^\circ$. Displayed is the average multigrid convergence rate of a V(1,1)-cycle with GMRES(3) smoother as a function of the energy $E$. Average convergence rate $(\|r_k\|/\|r_0\|)^{1/k}$ calculated from experimental results based upon $k = 4$ consecutive V-cycles. \label{fig:2dconvrate}}
\end{figure}

In this section, we benchmark the performance of multigrid as a solver
for the 2D Schr\"odinger scattering problem on a complex-valued
grid. It appears that the multigrid convergence rate critically
depends on the value of the total energy $E$.  Indeed, Figure
\ref{fig:2dconvrate} shows the convergence rate of a standard
multigrid V(1,1)-cycle for the 2D model problem described in the
sections above. We observe that the multigrid scheme fails to converge
for total energies $E$ between $-1$ and $0$. Note that this
corresponds precisely to the energy range where only single ionization
occurs. However, for energy levels $E>0$, where both single and double
ionization occur, multigrid successfully converges.

The observed convergence behavior can be explained using the spectral
properties of the Hamiltonian operator introduced in Section
\ref{subsec:spectralprops}. The bottom right panel of Figure
\ref{fig:eigenvalues} shows the eigenvalues of $H^{\text{1d}} \otimes
I + I \otimes H^{1d}$, which approximate the spectrum of
$H^{\text{2d}}$, discretized along a complex contour. Changing the
total energy $E$ shifts the spectrum of the linear system
$(H^{2d}-E)u=\phi$ to the left or right.  For $-1 < E < 0$ the spectrum
 shifts to the right, resulting in a spectrum with
real-valued eigenvalues both to the left and right of the origin.
This indeed implies difficulties for both the smoother and the coarse
grid correction scheme. However, when $0 < E$,
the spectrum shifts to the left, moving all eigenvalues away from
zero. This results in a spectrum that is distinctly separated from the
origin, corresponding to a problem more amenable to
iterative solution.

From figure \ref{fig:2dconvrate}, it might appear to the reader
that multigrid is not generally efficient as a Schr\"odinger solver
due to the poor convergence in the $-1 < E < 0$ region. However, it is
important to note that multigrid is performant in the region of
physical interest.  Indeed, the double ionization problem requires a
full two-dimensional description as stated above, requiring multigrid
to converge for energies $E>0$. In contrast, the purely single ionized
problem can be solved in a one-dimensional Helmholtz setting, see
\eqref{eq:boundstate}, where multigrid can indeed be shown to perform
well for energy levels $-1 < E < 0$, cf.\ Section \ref{sec:Numerical}.

Although single ionization waves are present in the solution for
energy levels $E > 0$, they do not undermine the multigrid convergence
in this regime. This is remarkable, because single ionization waves
are very localized evanescent waves along the edges of the domain that
generally cannot be represented efficiently on coarser grids, since
there might not be enough grid points covering these regions.
However, despite the fact that the coarsening strategy of the
multigrid method used in this work is not adapted to evanescent waves,
the damping implied by the complex contour evaluation ensures
good multigrid performance.

\subsection{Solutions of a 3D Schr\"odinger equation}\label{subsec:3dsystem}

As demonstrated on a 2D model problem in the previous sections, the far field map (cross section) of a general Schr\"odinger problem can be accurately
calculated using the complex contour approach. In this section we focus on the numerical
solution of the 3D Schr\"odinger equation on the complex contour using multigrid. 
Note that in the three-dimensional case, the use of a direct solver is strongly prohibited due to the size of the problem. 

We consider the 3D Schr\"odinger equation, modeling a
realistic scattering problem that includes single, double and triple
ionization. As discussed above, this problem features very
localized waves that require a sufficiently high-resolution representation.
The model problem derived from a 9D problem through a partial wave expansion is 
\begin{align}\label{eq:partialwave3d}
  \Big(-\frac{1}{2} \Delta  &+ V_1(x) +  V_2(y) + V_3(z) \notag \\
  &   + V_{12}(x,y) + V_{23}(y,z)+ V_{31}(z,x)-E \, \Big) \, u(x,y,z) = \phi(x,y,z),  \quad x,y,z \ge 0,
\end{align}
with boundary conditions
\begin{equation}\label{eq:partialwave3dBC}
\begin{cases} 
u(x,y,0)=0 \quad \text{for} \quad x,y\ge 0\\
u(0,y,z)=0 \quad \text{for} \quad y,z\ge 0\\
u(x,0,z)=0 \quad \text{for} \quad x,z\ge 0\\
\text{outgoing} \quad \text{for} \quad x \rightarrow \infty \quad \text{or} \quad y \rightarrow \infty \quad \text{or} \quad z \rightarrow \infty,
\end{cases}
\end{equation}
We discuss a system in which $V_1$, $V_2$ and $V_3$ are identical
one-body potentials and $V_{12}$, $V_{23}$ and $V_{31}$ are,
similarly, identical two-body potentials. Let the strength of the
one-body potential be such that there is a single negative eigenvalue
for the 1D subsystem
\begin{equation}\label{eq:1d_eigenstate}
\left(-\frac{1}{2} \frac{d^2}{dx^2}  + V(x)\right) \phi_0(x) = \lambda_0 \, \phi_0(x),   \quad x \ge 0,
\end{equation}
with $\lambda_0 <0$, where we have dropped the subscript on the 1D potential $V$.
If the two-body potential $V_{12}(x,y)$ is negligibly small, then there
automatically exists a bound state of the 2D subsystem. Indeed, the state
$\phi_0(x)\phi_0(y)$ is an eigenstate of the separable Hamiltonian
$(-1/2) \Delta + V(x) + V(y)$ with eigenvalue $2\lambda_0$. In the
presence of a small but non-negligible two-body potential, this state will be slightly
perturbed, resulting in an eigenstate $\phi_0(x,y)$ that fits the 2D subsystem
\begin{equation}\label{eq:2d_eigenstate}
\left(-\frac{1}{2} \Delta  + V(x) +  V(y) + V_{12}(x,y)\right) \phi_0(x,y) = \mu_0\phi_0(x,y),   \quad x,y \ge 0.
\end{equation}
The corresponding eigenvalue is $\mu_0 \approx 2\lambda_0
< \lambda_0 < 0$. This ordering is typical for realistic atomic and
molecular systems \cite{balslev1971spectral}.
Similarly, the 3D system will have an eigenstate that looks
approximately like $\phi_0(x,y)\phi_0(z)$, or any of its coordinate
permutations. This 3D eigenstate $\phi_0(x,y,z)$ fits the equation
\begin{align} \label{eq:3d_eigenstate}
  \Big( -\frac{1}{2} \Delta  &+ V_1(x) +  V_2(y) + V_3(z) \notag \\
  &   + V_{12}(x,y) + V_{23}(y,z)+ V_{31}(z,x) \Big) \, \phi_0(x,y,z) = \nu_0 \, \phi_0(x,y,z),  \quad x,y,z \ge 0,
\end{align}
where $\nu_0 \approx  \mu_0 + \lambda_0 \approx 3 \lambda_0$. 
 
\begin{table}
\begin{center}
 {\small
  \begin{tabular}{|c|c|c|c|c|c|c|c|c|c|}
\hline
     								   & ~~~$E < \nu_0$~~~~ & $\nu_0 < E < \mu_0$ & $\mu_0 < E < \lambda_0$ & $\lambda_0 < E < 0~$ & ~~~~$0 < E$~~~~  \\
\hline
     Indefinite        & 	No	& Yes & Yes & Yes & Yes \\
     Single ion. &	No	& No	& Yes & Yes & Yes \\
     Double ion. &	No	& No	&	No  & Yes & Yes \\
     Triple ion. &	No	& No  &	No	& No	& Yes \\
\hline
  \end{tabular}
  \vspace{0.2cm}
  }
\end{center}
\caption{
  Schematic overview of the different scattering regimes in the 3D model problem described in Section~\ref{subsec:3dsystem}.
  Depending on the value of the total energy $E$ with respect to the eigenvalues $\nu_0$, $\mu_0$ and $\lambda_0$ 
  of the 3D, 2D and 1D (sub-)system respectively, different types of scattering occur. 
  The system is indefinite as soon as $\nu_0< E$. Single ionization waves emerge as soon as $\mu_0 < E$ 
  and double ionization occurs for $\lambda_0 < E$. Triple ionization waves are present only when $0 < E$. 
}
\end{table}
  
Assuming that the potentials are such that $\nu_0 < \mu_0 < \lambda_0
< 0$, there are now four possible regimes of interest in equation
(\ref{eq:partialwave3d}), depending on the total energy $E$.  First,
for $E < \nu_0$, the problem is positive definite, and hence easy to
solve numerically. However, in this regime no interesting physical
reactions occur.  Similarly, for $ \nu_0 < E < \mu_0$, there are no
scattering states in the solution.  For energy levels $\mu_0 < E <
\lambda_0$, single ionization scattering occurs. Consequently, in this
regime, there exist scattering solutions that are localized along one
of the three axes in the 3D domain. These solutions take the form
$v(z)\phi_0(x,y)$ as $z \rightarrow \infty$, where $\phi_0(x,y)$ is
the eigenstate of \eqref{eq:2d_eigenstate} and $v(z)$ is a scattering
solution satisfying outgoing wave boundary conditions.  Similar
solutions are found for the respective coordinate permutations.  For
energies $\lambda_0 < E < 0$, both single and double ionization
occurs. The solution contains --- besides single ionization waves ---
double ionization waves of the form $w(y,z)\phi_0(x)$, where
$\phi_0(x)$ is an eigenstate of \eqref{eq:1d_eigenstate} and $w(y,z)$
is a 2D scattering state satisfying the outgoing wave boundary
conditions.  Together with coordinate permutations, these waves are
localized along the faces of the 3D domain, where one of the three
coordinates, $x$, $y$ or $z$, is small.  Finally, for $E>0$, the
solution contains, in addition, triple ionization waves. These are
waves that describe a quantum mechanical system that is fully broken
up into its sub-particles. In this case, all three relative
coordinates $x$, $y$ and $z$ can become large, resulting in a wave
that extends to the entire domain.

Note that in fact only the latter problem, when $E>0$ and triple
ionization is present, requires a full 3D description. For the regimes
in which only single ionization occurs, a description in the
form of a coupled set of 1D equations is sufficient, due to the
separated character of the solution.  Similarly, for problems with
both single and double ionization, but no triple ionization, a simpler
2D description such as the one given by \eqref{eq:partialwave2d} can
be used to fully describe the physics behind the problem. Hence, in
view of the efficient solution of the 3D Schr\"odinger problem
\eqref{eq:3d_eigenstate}, our main interest goes out to the $E > 0$
regime.

\subsection{Multigrid performance on the 3D Schr\"odinger problem} \label{subsec:MG3D}

We now study the convergence of a multigrid solver for the 3D Schr\"odinger 
equation (\ref{eq:partialwave3d}) for energies $E$ that cover all possible 
scattering regimes. The model problem under consideration is a straightforward 
generalization of the 2D model presented in Section \ref{subsec:solutypes}, featuring one-body 
potentials $V_1(x) = -4.5\exp(-x^2)$, $V_2(y) = -4.5\exp(-y^2)$ and 
$V_3(z) = -4.5\exp(-z^2)$, and two-body potentials 
$V_{12}(x,y) = 2\exp(-(x+y)^2)$, $V_{23}(y,z) = 2\exp(-(y+z)^2)$ and 
$V_{31}(x,z) = 2\exp(-(x+z)^2)$. These potentials imply the existence of a 1D 
eigenstate in \eqref{eq:1d_eigenstate} with corresponding energy $\lambda_0=-1.0215$, a 2D 
eigenstate solution of \eqref{eq:2d_eigenstate}
with energy $\mu_0=-1.841$, and an additional 3D eigenstate in
\eqref{eq:3d_eigenstate}, which has energy $\nu_0=-2.751$. The problem is 
solved for a range of different total energies $E$, using an identical
right-hand side $\phi(x,y,z) = \exp(-3(x+y + z)^2)$ for all energies.
The discretization comprises $255^3$ points, covering the complex-valued
cube domain $[0,15e^{i\pi/12}]^3$.

The 3D model problem described above is solved using a full multigrid F(5)-cycle
\cite{trottenberg2001multigrid}. This implies that the problem is first discretized
on a $7^3$-point grid, where it is solved exactly. The solution obtained on this level is consecutively
interpolated and used as an initial guess to the same problem discretized using $15^3$ grid points,
after which 5 V(1,1)-cycles are applied. This process is repeated recursively until
we arrive at the finest level in the multigrid hierarchy consisting of $255^3$ grid points.  
On this level the V-cycle convergence rate is measured by averaging the residual 
reduction rate over three consecutive V-cycles.
Figure \ref{fig:convergencerate} shows the convergence rate as a function of the total energy $E$. 
We observe acceptable convergence behavior for energy levels $E <-3.0$, where the problem is
positive definite.  However, for energy levels between $-3.0 < E < 0$, where single and double 
ionization occur, unacceptably slow convergence is measured.  For energy levels $E > 0$, where single, double and triple
ionization waves coexist, multigrid convergence is again good.

\begin{figure}
\begin{center}
\begin{tikzpicture}
\begin{axis}[
        xlabel=$E$,
        ylabel=Multigrid Convergence rate]
        \addplot[smooth,black] plot coordinates{
(   -4.0000   , 0.0689 )
(   -3.9000   , 0.0696 )
(   -3.8000   , 0.0704 )
(   -3.7000   , 0.0737 )
(   -3.6000   , 0.0809 )
(   -3.5000   , 0.0778 )
(   -3.4000   , 0.0741 )
(   -3.3000   , 0.0923 )
(   -3.2000   , 0.1086 )
(   -3.1000   , 0.1389 )
(   -3.0000   , 0.2434 )
(   -2.9000   , 0.8971 )
(   -2.8000   , 0.9972 )
(   -2.7000   , 1.0134 )
(   -2.6000   , 1.0067 )
(   -2.5000   , 0.9907 )
(   -2.4000   , 0.9746 )
(   -2.3000   , 0.9617 )
(   -2.2000   , 0.9592 )
(   -2.1000   , 0.9651 )
(   -2.0000   , 0.9783 )
(   -1.9000   , 0.9919 )
(   -1.8000   , 1.0014 )
(   -1.7000   , 1.0069 )
(   -1.6000   , 1.0103 )
(   -1.5000   , 1.0110 )
(   -1.4000   , 1.0109 )
(   -1.3000   , 1.0088 )
(   -1.2000   , 1.0041 )
(   -1.1000   , 0.9963 )
(   -1.0000   , 0.9849 )
(   -0.9000   , 0.9676 )
(   -0.8000   , 0.9380 )
(   -0.7000   , 0.8771 )
(   -0.6000   , 0.7819 )
(   -0.5000   , 0.6327 )
(   -0.4000   , 0.5303 )
(   -0.3000   , 0.4918 )
(   -0.2000   , 0.4629 )
(   -0.1000   , 0.4407 )
(         0   , 0.5800 )
(    0.1000   , 0.9002 )
(    0.2000   , 0.9485 )
(    0.3000   , 0.9542 )
(    0.4000   , 0.5505 )
(    0.5000   , 0.4018 )
(    0.6000   , 0.3580 )
(    0.7000   , 0.3462 )
(    0.8000   , 0.3349 )
(    0.9000   , 0.3216 )
(    1.0000   , 0.3177 )
(    1.1000   , 0.3096 )
(    1.2000   , 0.3020 )
(    1.3000   , 0.2951 )
(    1.4000   , 0.2881 )
(    1.5000   , 0.2810 )
(    1.6000   , 0.2746 )
(    1.7000   , 0.2691 )
(    1.8000   , 0.2643 )
(    1.9000   , 0.2599 )
(    2.0000   , 0.2559 )
(    2.1000   ,   0.2523    )
(    2.2000   ,    0.2488   )
(    2.3000   ,    0.2457   )
(     2.4000  ,    0.2428   )
(     2.5000  ,    0.2402   )
(     2.6000  ,    0.2379   )
(     2.7000  ,    0.2358   )
(     2.8000  ,    0.2340   )
(     2.9000  ,    0.2324   )
(     3.0000  ,    0.2311   )
(     3.1000  ,    0.2300   )
(     3.2000  ,    0.2292   )
(     3.3000  ,    0.2284   )
(     3.4000  ,    0.2279   )
(     3.5000  ,    0.2274   )
(     3.6000  ,    0.2270   )
(     3.7000  ,    0.2267   )
(     3.8000  ,    0.2264   )
(     3.9000  ,    0.2261   )
(     4.0000  ,    0.2258   )
(     4.1000  ,    0.2256   )
(     4.2000  ,    0.2253   )
(     4.3000  ,    0.2250   )
(     4.4000  ,    0.2247   )
(     4.5000  ,    0.2244   )
(     4.6000  ,    0.2241   )
(     4.7000  ,    0.2238   )
(     4.8000  ,    0.2234   )
(     4.9000  ,    0.2231   )
(     5.0000  ,    0.2227   )
(     5.1000  ,    0.2223   )
(     5.2000  ,    0.2220   )
(     5.3000  ,    0.2216   )
(     5.4000  ,    0.2212   )
(     5.5000  ,    0.2208   )
(     5.6000  ,    0.2204   )
(     5.7000  ,    0.2200   )
(     5.8000  ,    0.2195   )
(     5.9000  ,    0.2191   )
(     6.0000  ,    0.2186   )
(     6.1000  ,    0.2182   )
(     6.2000  ,    0.2177   )
(     6.3000  ,    0.2172   )
(     6.4000  ,    0.2167   )
(     6.5000  ,    0.2162   )
(     6.6000  ,    0.2157   )
(     6.7000  ,    0.2152   )
(     6.8000  ,    0.2146   )
(     6.9000  ,    0.2141   )
(     7.0000  ,    0.2135   )
(     7.1000  ,    0.2129   )
(     7.2000  ,    0.2123   )
(     7.3000  ,    0.2116   )
(     7.4000  ,    0.2109   )
(     7.5000  ,    0.2102   )
(     7.6000  ,    0.2094   )
(     7.7000  ,    0.2086   )
(     7.8000  ,    0.2077   )
(     7.9000  ,    0.2068   )
(     8.0000  ,    0.2058   )
};                         
\end{axis}
\end{tikzpicture}
\end{center}
\caption{
3D Schr\"odinger problem \eqref{eq:partialwave3d} for a total energy range $E \in [-4,8]$ solved on a full complex grid with $\gamma = \pi/12 = 15^\circ$. Displayed is the average multigrid convergence rate of a V(1,1)-cycle with GMRES(3) smoother as a function of the energy $E$. Average convergence rate $(\|r_k\|/\|r_0\|)^{1/k}$ calculated from experimental results based upon $k = 3$ consecutive V-cycles. 
\label{fig:convergencerate}}
\end{figure} 

In analogy to the 2D problem, the observed lack of convergence for a limited range of energies can be
understood by analyzing the spectral properties of the discrete 3D Schr\"odinger
operator.  The discrete operator comprises a Laplacian operator,
discretized along a complex contour, which results in multiple
series of eigenvalues stretching deep into the bottom half of the
complex plane. A more detailed discussion on the eigenvalues of the 
discrete Helmholtz operator can be found in \cite{reps2010indefinite}.
For the total Schr\"odinger operator, the eigenvalues differ slightly. 
For each subsystem eigenvalue, either 2D or 1D,
a series of eigenvalues arises from just below the
real axis into the negative half of the complex plane, cf.\ Figure \ref{fig:eigenvalues}. Changing the
energy $E$ shifts the distribution of these eigenvalues in the direction of the
real axis.  For energy levels $\nu_0 < E < 0$, there are series of eigenvalues
both in the third and the fourth quadrant in the complex plane.  Both
series start close to the real axis, resulting in an indefinite problem
with eigenvalues closely near the origin, causing poor multigrid convergence.
Contrarily, in the $0 < E$ regime, all eigenvalues are bounded away from
the origin. Indeed, in this case all eigenvalue series start from eigenvalues along
the negative real axis. The largest real-valued eigenvalue lies at a distance
$|E|$ to the left of origin, implying the entire spectrum can be distinctly separated from 
the origin by a virtual straight line. This spectral property results in good multigrid convergence. 

Note, however, that from a physical point of view, the lack of convergence for energy levels 
$E < 0$ is not a concern, since the solutions in this energy regime
can be described either by a 1D or 2D equation (see higher), and a full 3D
description is generally not required.

%
%
%
%

\section{Conclusions and discussion}\label{sec:Conclusions}
In this paper we have developed a novel highly efficient method for
the calculation of the far field map resulting from $d$-dimensional
Helmholtz and Schr\"o\-ding\-er type scattering problems where the
wavenumber is an analytical function. Our approach is based on the
reformulation of the classically real-valued Green's function volume
integral for the far field map to an equivalent volume integral over a
complex-valued domain.

The advantage of the proposed reformulation lies in the scattered wave
solution of the Helmholtz problem on a complex domain, which can be
calculated efficiently using a multigrid method. This is particularly advantageous
for 3D problems, where direct solution is generally very hard.  Indeed, the
reformulation of the Helmholtz forward problem on the full complex
contour is shown to be equivalent to a Complex Shifted Laplacian
problem, for which multigrid has been proven in the literature to be a
fast and scalable solver. However, whereas the Complex Shifted
Laplacian was previously only used as a preconditioner, the
complex-valued far field map calculation proposed in this paper
effectively allows for multigrid to be used as a solver on the
perturbed problem.

The functionality of the method is primarily validated on 2D and 3D
Helmholtz type model problems. It is confirmed that the values of the
far field map calculated on the full complex grid exactly matches the
values of the classical real-valued integral. Furthermore, the number
of multigrid iterations is shown to be largely wavenumber independent,
yielding a fast overall far field map calculation.

We have found that rotating the contour about 15 degrees into the
complex domain is sufficient to ensure multigrid stability. Choosing a 
larger rotation angle improves the multigrid convergence, but makes 
the far field integral harder to calculate since it is a product of an
exponentially decaying with a exponentially increasing function, and a
larger angle implies an increased rate of decay or growth. The
limited machine precision of the implementation may influence the accuracy 
of the integrand product for large domains and angles.

One area of scientific computing where the proposed technique might be
particularly valuable is in the numerical solution of quantum
mechanical scattering problems. These are generally high-dimensional
scattering problems where the wavenumber is indeed an analytical
function, and where 6D or 9D problems are common. We have validated
that for a 2D Schr\"odinger type model problem the proposed method is able to
accurately calculate the cross sections that are measured in physical
experiments. In addition, we have studied the convergence rate for a
3D Schr\"odinger equation. Results show reasonable multigrid
convergence rates for the energy range of interest.

Despite their use as benchmark problems, the model problems described
in this paper form an important testing framework for more realistic
applications.  However, further analysis of the convergence rates are
necessary for realistic Coulomb potentials to make the method more
robust, and eventually usable by computational physicists and
chemists.

Note that the proposed method can also be used to calculate near-field
amplitudes that play an important role in lithography, microscopy and
tomography.  These amplitudes are calculated in the same way as the
far field map as an integral over the Green's function and the
numerical solution. For the near-field however, no asymptotic form of the
Green's function can be used. In principle, this integral can also be 
calculated along a complex valued contour and will yield convergence 
results similar to the far field results presented in this work.

Finally, we note that a number of modifications can be made to improve
the efficiency of the method even further, like choosing the shape of
the complex contour for the integral based on a steepest descent
scheme, as proposed in \cite{huybrechs2006evaluation}.

\section{Acknowledgments}
This research was partly funded by the \textit{Fonds voor
  Wetenschappelijk Onderzoek (FWO)} project G.0.120.08 and
\textit{Krediet aan navorser} project number 1.5.145.10. Additionally,
this work was partly funded by
Intel$^{{\scriptsize\textregistered}}\hspace{-0.05cm}$ and by the
\textit{Institute for the Promotion of Innovation through Science and
  Technology in Flanders (IWT)}. The authors would like to thank
Hisham bin Zubair for sharing a multigrid implementation and D.
Huybrechs, C.W.\ McCurdy and D.J.\ Haxton for fruitful discussions on
the subject.

\nocite{*}
{\small
\bibliographystyle{plain}
\bibliography{refs}
}

\end{document}